\documentclass[aos,11pt,noinfoline]{imsart}
\usepackage{booktabs}
\usepackage[figuresright]{rotating}
\usepackage{fancybox}
\usepackage{curves}
\usepackage{latexsym}
\usepackage{amssymb}

\usepackage{graphicx}
\usepackage{amsfonts}
\RequirePackage{amsthm,amsmath}
\RequirePackage[numbers,sort&compress]{natbib}
\RequirePackage[colorlinks,citecolor=blue,urlcolor=blue]{hyperref}
\usepackage{epic}
\usepackage{eepic}

\startlocaldefs

\newcommand{\newmathbin}[2]{\def#1{\mathbin{#2}}}

\newcommand{\E}[1]{{\mathop{\mathrm{E}\kern-1pt}\nolimits}\left(#1\right)}
\newcommand{\Cov}[1]{{\mathop{\mathrm{Cov}\kern-1pt}\nolimits}\left(#1\right)}
\newcommand{\Var}[1]{{\mathop{\mathrm{Var}\kern-1pt}\nolimits}\left(#1\right)}
\newcommand{\trace}{\mathop{\mathrm{trace}}}
\newcommand{\rank}{\mathop{\mathrm{rank}}}
\newcommand{\image}{\mathop{\mathrm{Im}}\nolimits}
\newcommand{\setof}[1]{\ensuremath{\left\{#1\right\}}}
\newcommand{\setofall}[2]{\setof{#1:#2}}
\newcommand{\card}[1]{\left|#1\right|}
\newcommand{\combine}{\mathbin{\vartriangleright}}

\newcommand{\resid}{\mathbin{\vdash}}
\newcommand{\sresid}{\mathbin{\vdash}}
\newmathbin{\joint}{\square}
\newcommand{\inter}{\mathbin{\#}}
\newcommand{\nesting}[1]{\left[#1\right]}
\newcommand{\meet}{\wedge}

\newcommand{\summ}{\mathop{\sum\nolimits'}}

\newcommand{\add}{\mathbin{+}}

\newcommand{\modsep}{\mathbin{\mid}}

\newcommand{\blob}{\circle*{0.2}}

\newcommand{\nonorthcircle}{\put(0,0){\circle{0.2}}}
\newcommand{\orthcircle}{\put(0,0){\circle{0.2}}\put(-0.11,-0.06){{\tiny$\bot$}}}

\newcommand{\makepseudo}{\put(0,0){$\blacklozenge$}}

\newlength{\nlevnamesep}
\newsavebox{\savetier}
\newenvironment{tierbox}{\begin{lrbox}{\savetier}
\begin{tabular}{r@{\hspace{\nlevnamesep}}l}}{\end{tabular}\end{lrbox}
\ovalbox{\usebox{\savetier}}}

\theoremstyle{definition}
\newtheorem{egg}{Example}

\newenvironment{enumeroman}{
\begin{enumerate}}
{\end{enumerate}}

\endlocaldefs

\oddsidemargin
\heavyrulewidth=\lightrulewidth

\begin{document}

\begin{frontmatter}

\title{Randomization-based models for multitiered experiments:
I. A chain of randomizations}
\runtitle{Randomization-based models}

\begin{aug}
\author{\fnms{R. A.} \snm{Bailey}\thanksref{m1}\corref{}\ead[label=e1]
{rab@mcs.st-and.ac.uk}}
\and
\author{\fnms{C. J.} \snm{Brien}\thanksref{m2}\ead[label=e2]{chris.brien@unisa.edu.au}
\ead[label=u1,url]{http://chris@brien.name}}

\address{Professor R. A. Bailey \\
         School of Mathematics and Statistics\\
         University of St Andrews\\
         St Andrews,
         Fife KY16 9SS\\
         United Kingdom \\
\printead{e1}\\ \\
         School of Mathematical Sciences \\
         Queen Mary University of London \\
         Mile End Road \\
         London  E1 4NS \\
         United Kingdom}

\address{Dr C. J. Brien \\Phenomics and Bioinformatics Research Centre \\
School of Information Technology and Mathematical Sciences \\
University of South Australia \\
GPO Box 2471, Adelaide, SA 5001 \\ \printead{e2}\\ \\
The Australian Centre for Plant Functional Genomics \\
Waite Campus, University of Adelaide \\
Urrbrae, SA 5064\\
Australia \\\printead{u1}}

\affiliation{University of St Andrews\thanksmark{m1} and 
Queen Mary, University of London\thanksref{m1}, 
University of South Australia\thanksmark{m2} and 
University of Adelaide\thanksmark{m2}}

\runauthor{R.A. Bailey and C.J. Brien}
\end{aug}

\today

\begin{abstract}
We derive randomization-based models for experiments with
a chain of randomizations.
The estimation theory for these models leads to formulae for the
estimators of treatment effects, their standard errors, and
expected mean squares in the analysis of variance.
We discuss the practicalities in
fitting these models and outline the difficulties that can occur,
many of which do not arise in two-tiered experiments.
\end{abstract}



\end{frontmatter}

\section{Introduction}

\citet{Bailey81,Bailey91}, following \citet{Grundy50}, outlines a
method of deriving randomization-based models for
experiments. 
It provides mixed models that are randomization-based in the sense that
it is the group of permutations for the randomization that defines the
variance matrix on which the analysis is to be based.
It applies to a general class of structures: those derived from a group
of permutations. Although this approach is
restricted to groups which are stratifiable, in the sense defined in
Section~\ref{sec:ranmod}, it includes all poset block structures and many
other structures besides.

A restriction with this approach is that it only applies to a
single randomization,
as defined by \citet{BrBa:mult}, in that the randomization can be
achieved using a single permutation of the set of observational units.
Brien and Bailey \citep{BrBa:mult} 
describe experiments with multiple randomizations, requiring multiple
permutations, and show how to assess 
the properties of such experiments in \citep{Brien092,Brien10}.
As for the analysis of such experiments,
Curnow \citep{Curnow59}, in correcting the analysis of
McIntyre \citep{McIntyre55}, showed how to analyse the results of
two-phase experiments by analysis of variance (anova).
Wood, Williams and Speed \citep{Wood88} also discussed
the analysis of two-phase experiments.
Brien in \citep{Brien83} indicated how to
use tiers to obtain the anova for multitiered experiments
and in \citep{Brien92s} derived
expected mean squares under a mixed model.
Brien and Payne \citep{Brien99} extended the sweep algorithm of Wilkinson
\citep{Payne77,Wilkinson70} to cover anova for multitiered experiments.
Brien and Bailey \citep{BrBa:mult} and Brien and Dem\'etrio \citep{Brien091}
describe how to analyse the data from such experiments by using mixed models.
However, no one has so far given general formulae for the estimators
of treatment effects and their standard errors for multitiered
experiments, nor have formulae for the expected mean squares
under randomization-based models been derived.

Section~\ref{sec:ranmod} formulates the randomization-based model
for a two-tiered experiment and 
generalizes it to experiments with two randomizations in a chain.
Section~\ref{sec:treatsb} describes families of expectation models that
lead to a treatment decomposition;
the assumption of structure balance is discussed.
The properties of anovas corresponding to randomization-based models
are outlined in Section~\ref{sec:sameaov}. Section~\ref{sec:eg}
contains a set of examples. 
Sections~\ref{sec:esti2}--\ref{sec:comk} address the
estimation of treatment effects, first for two-tiered experiments and then for
various cases of three-tiered experiments.
Section~\ref{sec:multi} generalizes this to an arbitrary number of
randomizations in a chain. Section~\ref{sec:dataest}
covers the use of software in estimating model parameters, 
including a discussion of randomization-based models in the class
of all mixed models.
Statistical inference is discussed in Section~\ref{sec:infer}.
Section~\ref{sec:altmod} briefly touches on models other than those described
in Section~\ref{sec:ranmod}. See Section~3 of \citep{Brien092} for
definitions of some terms and notation specific to the approach we take.

\section{Randomization-based models}
\label{sec:ranmod}

\subsection{The randomization-based model for a two-tiered experiment}
\label{sec:rand2}

As in \citep{BrBa:mult,Brien092,Brien10},
we randomize the set of
objects~$\Gamma$ to another set of objects~$\Upsilon$, so we have
a design function $h\colon\Upsilon\to\Gamma$.  If the objects in
$\Gamma$ are treatments then $h(\upsilon)$ is the treatment
assigned to unit $\upsilon$ in $\Upsilon$.
We associate a structure with each of $\Upsilon$ and $\Gamma$.
If $V_\Upsilon$ is the space of all real vectors indexed by
$\Upsilon$, then a \emph{structure} on $\Upsilon$ is an orthogonal
decomposition of $V_{\Upsilon}$.  This is specified by
a set of  symmetric, idempotent, mutually orthogonal matrices
projecting onto the subspaces of~$V_{\Upsilon}$ in the
decomposition. Similarly, structure on $\Gamma$ is an orthogonal
decomposition of the space $V_\Gamma$.

The usual initial assumption for the response~$Y_\upsilon$ on
unit~$\upsilon$ in~$\Upsilon$ is additive:
\begin{equation}
  Y_\upsilon = w_\upsilon + \tau_{h(\upsilon)}.
\label{eq:add}
\end{equation}
In some approaches, 
$w_\upsilon$ is taken to be a constant, but
here it is taken to be  a random variable, as in
\citep{Bailey91,Grundy50}. 
It depends only on the unit $\upsilon$ which is
providing the response.  On the other hand, $\tau_i$, for $i$ in
$\Gamma$, is a constant. It depends only on the treatment~$i$
which is applied to $\upsilon$. 
Permitting the $w_\upsilon$ to be random allows for
measurement error,  without the assumption of any particular form
for it, and any random sampling of units that may occur.

Let $G$ be a group of permutations of $\Upsilon$.  We usually take $G$ to
be the largest group of permutations that preserve certain
generalized factors on $\Upsilon$, in the sense that if $F$ is
such a generalized factor and $g\in G$ and $F(\upsilon_1) = F(\upsilon_2)$
then we must have $F(g(\upsilon_1)) = F(g(\upsilon_2))$.
In \citep{Bailey81,Bailey91} it is argued that if we randomize by choosing
$g$ from $G$ at random then it is appropriate to replace $w_\upsilon$ by
$W_\upsilon$, which is the mixture of the $w_{g(\upsilon)}$ over $g$ in
$G$.  Hence we get
the randomization-based model
\begin{equation}
  Y_\upsilon = W_\upsilon + \tau_{h(\upsilon)},
\label{eq:rand}
\end{equation}
  where the $W_\upsilon$ are random variables
which are exchangeable under~$G$: in particular,
\begin{enumerate}
    \item[(P.a)]
if there is any $g$ in $G$ for which $g(\upsilon_1) = \upsilon_2$ then
$W_{\upsilon_1}$ and  $W_{\upsilon_2}$ have the same distribution, in
particular, the same expectation;
\item[(P.b)]
if there is any $g$ in $G$ for which $g(\upsilon_1) = \upsilon_2$
and $g(\upsilon_3) = \upsilon_4$ then the joint distribution of
$(W_{\upsilon_1},W_{\upsilon_3})$ is the same as the joint distribution
of $(W_{\upsilon_2},W_{\upsilon_4})$, in particular,
$\Cov{W_{\upsilon_1},W_{\upsilon_3}} = \Cov{W_{\upsilon_2},W_{\upsilon_4}}$.
\end{enumerate}
If the group $G$ is transitive on $\Upsilon$ then
property~(P.a) is true for all choices of $\upsilon_1$ and $\upsilon_2$, so
we may incorporate the constant value of $\E{W_\upsilon}$ into each
$\tau_i$ and so assume that $\E{W_\upsilon} = 0$ for all $\upsilon$ in
$\Upsilon$.  We restrict attention to cases where $G$ is transitive,
which implies that every unrandomized factor on $\Upsilon$ is equireplicate.

Let $\mathbf{Y}$ and $\mathbf{W}$ be the vectors of the random variables
$Y_\upsilon$ and $W_\upsilon$ respectively.
We often represent the design function $h$ by a design matrix
$\mathbf{X}_h$.  This is an $\Upsilon \times \Gamma$ matrix with
$(\upsilon,i)$-entry equal to $1$ if $h(\upsilon)=i$ and to $0$ otherwise.
Then equation~\eqref{eq:rand} can be rewritten in vector form as
$
\mathbf{Y} = \mathbf{W} + \mathbf{X}_h \boldsymbol{\tau}
$,
and
$\E{\mathbf{Y}} = \mathbf{X}_h \boldsymbol{\tau}$.

The pattern in the (co)variance matrix $\mathbf{C}$ 
of $\mathbf{W}$ is determined by property~(P.b), which implies that 
$\mathbf{C}$ is a patterned matrix with the same entries, including 
multi\-plicities, in every row; only their order differs. 
More specifically, there is a set $\mathcal{B}$
of symmetric $\Upsilon \times \Upsilon$ adjacency matrices $\mathbf{B}$
with entries $0$ and $1$, whose sum is the all-$1$ matrix $\mathbf{J}$,
such that if the $(\upsilon_1,\upsilon_2)$-entry of $\mathbf{B}$
is equal to $1$ then the $(\upsilon_3,\upsilon_4)$-entry is
equal to $1$ if and only if there is some $g$ in $G$ for which
either $g(\upsilon_1) = \upsilon_3$ and $g(\upsilon_2) = \upsilon_4$
or $g(\upsilon_1) =\upsilon_4$ and $g(\upsilon_2) =\upsilon_3$.
Moreover, the product of any two adjacency matrices is a linear
combination of matrices in $\mathcal{B}$.  Property~(P.b) implies that 
there are (co)variances $\zeta_{\mathbf{B}}$ such that
$
\mathbf{C}  = \sum_{\mathbf{B}\in\mathcal{B}} \zeta_{\mathbf{B}}\mathbf{B}
$.
For simple orthogonal block structures, this form of $\mathbf{C}$ 
is the same as the variance matrix for
the null randomization distribution given by \citet{Nelder65a}.

The group $G$ is said to be \emph{stratifiable} \citep{Alejandro03,BarAz:rand}
if the eigenvectors of the matrix~$\mathbf{C}$
do not depend on the values of the entries $\zeta_{\mathbf{B}}$
but depend only on their pattern.
Then the common eigenspaces of $\mathbf{C}$,
called \emph{strata}, are the structure on $\Upsilon$, and the 
collection of possible variance matrices is said
to have \emph{orthogonal variance structure} (OVS).
OVS is called `orthogonal block structure'
by Houtman and Speed in \citep{Houtman83}.
Note that there is no linear dependence among the
(co)variances $\zeta_{\mathbf{B}}$.
Unless otherwise stated, we assume that $G$ is stratifiable
and so $\mathbf{C}$ has OVS.
Then the number of strata is equal to the number of adjacency matrices.

Let $\mathcal{Q}$ be the collection of symmetric,
mutually orthogonal, idempotent matrices projecting onto the strata.
Then $\sum_{\mathbf{Q}\in\mathcal{Q}} \mathbf{Q}$  is the 
$\Upsilon\times\Upsilon$ identity matrix
$\mathbf{I}_{\Upsilon}$, and the variance matrix can be re-expressed as
\begin{equation} 
\mathbf{C} = \sum_{\mathbf{Q}\in\mathcal{Q}} \eta_{\mathbf{Q}} \mathbf{Q},
\label{eq:cov}
\end{equation}
with $\eta_{\mathbf{Q}}\geq 0$ for all $\mathbf{Q}$ in $\mathcal{Q}$.
The values $\eta_{\mathbf{Q}}$ are the eigenvalues of
$\mathbf{C}$ and are called \emph{spectral components of variance}. 
The strata are subspaces within which all normalized
contrasts have equal variance under randomization, this variance being
the $\eta$ for that stratum. Given $\mathcal{Q}$, any two  matrices of the
form~\eqref{eq:cov} commute with each other.
The matrices $\mathbf{Q}$ are linear combinations of the matrices
$\mathbf{B}$, and vice-versa, but in general there is no closed-form
expression for the coefficients in these combinations.

\subsection{Application to poset block structures}
\label{sec:pbs2}
The majority of experiments conducted in practice, and all the 
subsequent examples in this paper, have poset block structures
on their units.
All poset block structures have stratifiable permutation groups,
as shown in \citep{BaileyPRS83}.
A poset block structure on $\Upsilon$ is defined by
a set $\mathcal{H}$ of generalized factors on $\Upsilon$: see \citep{BrBa:mult}.
Following \cite{Tjur84d}, we write $H < F$ if $H$ and $F$ are
in $\mathcal{H}$ and $H$ is marginal to $F$.
There are several ways in which to write $\mathbf{C}$,
in terms of matrices and coefficients that depend on $H$ in
$\mathcal{H}$ \citep{Speed87a,Speed87}:
\begin{equation}
\mathbf{C} = \sum_{H\in\mathcal{H}} \zeta_H \mathbf{B}_H
           = \sum_{H\in\mathcal{H}} \psi_H \mathbf{S}_H
           = \sum_{H\in\mathcal{H}} \eta_H \mathbf{Q}_H.
\label{eq:covfact}
\end{equation}
Here $\zeta_H$ is the (co)variance under the randomization between
elements of $\Upsilon$ with
the same level of $H$ but not the same level of
any generalized factor $F$ in $\mathcal{H}$ to which $H$ is marginal;
$\mathbf{B}_H$ is the $\Upsilon\times\Upsilon$
adjacency matrix with entry~$1$ for such pairs
 and entry~$0$ otherwise;
$\psi_H$~is a canonical component;
$\mathbf{S}_H$ is the
$\Upsilon\times\Upsilon$ relationship matrix \citep{James57} for $H$,
with $(\upsilon_1, \upsilon_2)$-entry equal to~$1$ if $\upsilon_1$
and $\upsilon_2$ have the same level of $H$ and to $0$ otherwise.
Thus $\mathbf{S}_H = \sum_{F\geq H}\mathbf{B}_{F}$.
If $k_H$ is the common replication of all levels of $H$, then
$k_H^{-1}\mathbf{S}_H = \sum_{F\leq H}\mathbf{Q}_{F}$.
When $H$ is the generalized factor consisting of all factors on
$\Upsilon$, the subscript $H$ will sometimes be replaced by $\Upsilon$,
while the subscript for the generalized factor corresponding to the 
overall mean is denoted $0$.
Expressions in \citep{Speed87a,Speed87} show how to convert one set of
coefficients in equation~\eqref{eq:covfact} to another. In particular,
\begin{equation}
\eta_H = \sum_{F\geq H}  k_{F} \psi_{F}.
\label{eq:zetaphi}
\end{equation}

The natural interpretation of \emph{canonical components} in
this context is as components of excess covariance \citep{Nelder77}.
They are linear combinations of the covariances ($\zeta$ parameters)
\citep{Nelder65a,Nelder77,Speed87}.
Except for $\psi_\Upsilon$, which is the variance of each individual response,
they measure the covariation, between the
responses on the units in $\Upsilon$, contributed by each particular 
generalized factor in excess of that of any 
generalized factor which is marginal to it.
Thus $\psi_H$ can negative when $\zeta_H <\zeta_F$ and  $F<H$,
although $\psi_{\Upsilon} > 0$. This is in contrast to nonnegative
variance components $\sigma_H^2$, 
which occur in the usual formulations of mixed models
that we discuss in Section~\ref{sec:dataest}.
Estimates of standard errors of treatment effects require estimates of
the spectral components. On the other hand, scientifically interesting
hypotheses about the canonical components  are often formulated and
tested \citep{Cox58b,Nelder77} (see also Section~\ref{sec:mixmod}) 
and so estimates of them may also be required.

\subsection{The randomization-based
model for an  experiment with two randomizations in a chain}
\label{sec:derive}

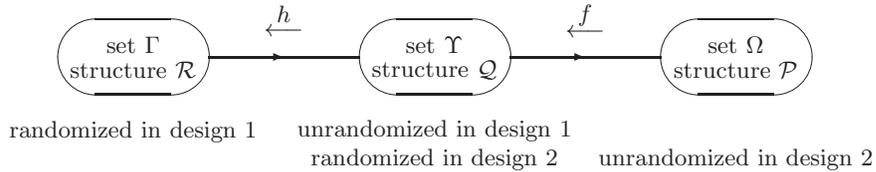
\begin{figure}[htbp]
\centering
\footnotesize
\setlength{\unitlength}{0.95pt}
\begin{picture}(350,70)(-175,-40)
\multiput(-120,5)(120,0){3}{\oval(60,30)}
\put(120,10){\makebox(0,0){set $\Omega$}}
\put(120,0){\makebox(0,0){structure $\mathcal{P}$}}
\put(-120,10){\makebox(0,0){set $\Gamma$}}
\put(-120,0){\makebox(0,0){structure $\mathcal{R}$}}
\put(0,10){\makebox(0,0){set $\Upsilon$}}
\put(0,0){\makebox(0,0){structure $\mathcal{Q}$}}
\put(30,5){\line(1,0){60}}
\put(-90,5){\vector(1,0){30}}
\put(-30,5){\line(-1,0){60}}
\put(30,5){\vector(1,0){30}}
\put(-60,20){\makebox(0,0){$\stackrel{\textstyle h}{\longleftarrow}$}}
\put(60,20){\makebox(0,0){$\stackrel{\textstyle f}{\longleftarrow}$}}
\put(0,-30){\makebox(0,0){\shortstack{unrandomized in design 1\\
randomized in design 2}}}
\put(120,-34){\makebox(0,0){\shortstack{ \\ unrandomized in design 2}}}
\put(-120,-26){\makebox(0,0){\shortstack{randomized in design 1\\ \hspace{1pt} }}}
\end{picture}
\caption{Diagram of an experiment with two randomizations in a chain}
\label{fig:samedir}
\end{figure}

For a chain of two randomizations, there are three sets: 
$\Upsilon$ is randomized to $\Omega$,
and $\Gamma$ is randomized to $\Upsilon$.  Let the corresponding
design maps be $f\colon\Omega\to\Upsilon$ and
$h\colon\Upsilon\to\Gamma$, as in Figure~\ref{fig:samedir}.
The elements of $\Gamma$ will be referred to as 
treatments and $\Upsilon$ and $\Omega$ as unrandomized sets.

Suppose that $f$~is randomized by choosing a random permutation
from the group~$G_1$ of permutations of~$\Omega$, and that $G_1$~is
stratifiable with stratum projectors~$\mathbf{P}$, for $\mathbf{P}$
in~$\mathcal{P}$.   Like the matrices
$\mathbf{Q}$ in Section~\ref{sec:rand2}, the matrices
$\mathbf{P}$ are known orthogonal idempotents summing to
the $\Omega \times \Omega$ identity matrix $\mathbf{I}_{\Omega}$.
In Section~\ref{sec:rand2} the size of the idempotents is the size of
$\Upsilon$, while here  it is the size of $\Omega$.

Now let $Y_\omega$ be the response on observational unit~$\omega$
in $\Omega$. 
Applying the randomization argument from Section~\ref{sec:rand2} to $f$ gives
$
Y_\omega = Z_\omega + \tilde Y_{f(\omega)}
$,
where $Z_\omega$ is a random variable depending only on the 
unit $\omega$ and $\tilde Y_\upsilon$ is a notional effect associated
with unit~$\upsilon$ in $\Upsilon$. Because $G_1$ 
is stratifiable, we can assume that the random
variables $Z_\omega$ are identically distributed with mean zero,
and that $\Cov{\mathbf{Z}} = \sum_{\mathbf{A}\in\mathcal{A}}
\gamma_{\mathbf{A}} \mathbf{A}$, where $\mathcal{A}$ is the set of
adjacency matrices arising from $G_1$ and the $\gamma_{\mathbf{A}}$
are the associated (co)variances. Following Section~\ref{sec:rand2},
we can also write $\Cov{\mathbf{Z}} = \sum_{\mathbf{P}\in\mathcal{P}}
\xi_\mathbf{P} \mathbf{P}$ where, like the quantities 
$\eta_{\mathbf{Q}}$, the $\xi_\mathbf{P}$ are unknown nonnegative coefficients.
Then $\Cov{\mathbf{Z}}$ has OVS because $G_1$ is stratifiable.

Similarly, $h$~is randomized by choosing a random permutation from the 
group~$G_2$ of permutations of~$\Upsilon$, and
$G_2$~is stratifiable with
$\Upsilon \times \Upsilon$ stratum projectors~$\mathbf{Q}$, for
$\mathbf{Q}$ in~$\mathcal{Q}$, as in Section~\ref{sec:rand2}.
Rewriting equation~\eqref{eq:rand} as
$
\tilde  Y_\upsilon = W_\upsilon + \tau_{h(\upsilon)}
$
gives 
\begin{equation}
Y_\omega = Z_\omega + \tilde Y_{f(\omega)}
         = Z_\omega + W_{f(\omega)} + \tau_{h(f(\omega))}.
\label{eq:2rand}
\end{equation}
In turn, this randomization-based model can be rewritten in vector form as
$
\mathbf{Y} = \mathbf{Z} + \mathbf{X}_f\mathbf{W} + \mathbf{X}_f
\mathbf{X}_h \boldsymbol{\tau}
$,
where $\mathbf{X}_f$ is the $\Omega \times \Upsilon$ design matrix for~$f$.
Hence $\E{\mathbf{Y}} = \mathbf{X}_f\mathbf{X}_h \boldsymbol{\tau}$,
and the variance matrix $\mathbf{V}$ of $\mathbf{Y}$ is given by
\[
\mathbf{V} = \Cov{\mathbf{Z} + \mathbf{X}_f\mathbf{W}}
           = \sum_{\mathbf{A}\in\mathcal{A}} \gamma_{\mathbf{A}} \mathbf{A} +
              \sum_{\mathbf{B}\in\mathcal{B}} \zeta_{\mathbf{B}}
                                      \mathbf{X}_f \mathbf{B} \mathbf{X}'_f
= \sum_{\mathbf{P}\in\mathcal{P}} \xi_\mathbf{P} \mathbf{P} +
\sum_{\mathbf{Q}\in\mathcal{Q}} \eta_\mathbf{Q} \mathbf{X}_f
\mathbf{Q} \mathbf{X}_f',
\]
because $\mathbf{Z}$ and $\mathbf{W}$ are independent.
The two sets $\mathcal{P}$ and $\mathcal{Q}$ of idempotents 
correspond to the eigenspaces of the variance matrices of $\mathbf{Z}$ 
and $\mathbf{W}$, respectively, but not necessarily to those of $\mathbf{V}$. 
Although the coefficients $\xi_{\mathbf{P}}$ and $\eta_{\mathbf{Q}}$ may not be
eigenvalues of $\mathbf{V}$, we still call them spectral components of
variance because they are the eigenvalues of
$\Cov{\mathbf{Z}}$ and $\Cov{\mathbf{W}}$ respectively.

As noted in \citep{Brien092},
the effect of the design function $f$ is to embed a copy $V_\Upsilon^f$ of
$V_\Upsilon$ inside  the space $V_\Omega$ of real vectors indexed by~$\Omega$.
Let $\mathbf{D}_f$ be the $\Upsilon \times \Upsilon$ diagonal matrix
whose $(\nu,\nu)$-entry is the replication of unit~$\nu$. Then
$\mathbf{X}'_f \mathbf{X}_f = \mathbf{D}_f$,
and the matrix of orthogonal projection onto $V_\Upsilon^f$ is
$\mathbf{X}_f \mathbf{D}_f^{-1} \mathbf{X}'_f$.

To further simplify $\mathbf{V}$,
the design $f$~be must be equireplicate.
If $\mathbf{Q}_1$ and $\mathbf{Q}_2$ are in $\mathcal{Q}$ then
$(\mathbf{X}_f \mathbf{Q}_1 \mathbf{X}_f')
(\mathbf{X}_f \mathbf{Q}_2 \mathbf{X}_f') =
\mathbf{X}_f \mathbf{Q}_1 \mathbf{D}_f \mathbf{Q}_2 \mathbf{X}_f'$.
If the common replication
is~$r$ then $\mathbf{D}_f = r\mathbf{I}_\Upsilon$, so if we put
$\mathbf{Q}^f= r^{-1} \mathbf{X}_f \mathbf{Q} \mathbf{X}_f'$ then
the $\mathbf{Q}^f$ are mutually orthogonal idempotents summing to
$r^{-1}\mathbf{X}_f\mathbf{X}_f'$, which is the matrix of orthogonal
projection onto the subspace $V_\Upsilon^f$.
To simplify notation, as in \citep{Brien092} we shall write 
$\mathbf{Q}^f$ just as $\mathbf{Q}$,
$\setofall{\mathbf{Q}^f}{\mathbf{Q}\in\mathcal{Q}}$ as $\mathcal{Q}$,
and $r^{-1}\mathbf{X}_f\mathbf{X}_f'$ as $\mathbf{I}_\mathcal{Q}$
in the three-tiered context.  Thus we have
\begin{equation}
\mathbf{V} = 
\sum_{\mathbf{P}\in\mathcal{P}} \xi_\mathbf{P} \mathbf{P} + r
\sum_{\mathbf{Q}\in\mathcal{Q}} \eta_\mathbf{Q} \mathbf{Q}.
\label{eq:twocov}
\end{equation}

The formula for $\mathbf{V}$ in equation~\eqref{eq:twocov} appears very
similar to that in equation~(2) of \citep{Wood88}.
There are three differences. Here, the two collections of idempotents 
sum to $\mathbf{I}$ and $\mathbf{I}_\mathcal{Q}$ respectively, whereas
those in \citep{Wood88} both sum to $\mathbf{I}$.
Equation~\eqref{eq:twocov} is justified by the randomization;
the formula in \citep{Wood88} is an assumed model.
Moreover, \citep{Wood88} does not require $f$ to be equireplicate, so
the parametrization does not explicitly include the replication~$r$.

\subsection{Pairs of poset block structures}
\label{sec:pbs3}

If the structure on $\Upsilon$ is a poset block structure with set 
$\mathcal{H}_2$ of generalized factors,
then $r\sum_{H\in\mathcal{H}_2} \eta_H \mathbf{Q}_H^f =
\sum_{H\in\mathcal{H}_2} \psi_H \mathbf{S}_H^f$, where
$ \mathbf{S}_H^f = \mathbf{X}_f  \mathbf{S}_H \mathbf{X}_f'$
which is the $\Omega \times \Omega$ relationship matrix for $H$ when
it is regarded as a factor on $\Omega$, in which case the common
replication of its levels is $rk_H$.  If we now write $\mathbf{S}_H^f$
just as $\mathbf{S}_H$, we have 
$r \sum_{H\in\mathcal{H}_2} \eta_H \mathbf{Q}_H 
= \sum_{H\in\mathcal{H}_2} \psi_H \mathbf{S}_H$.

Suppose that the structure on $\Omega$ is also a poset block
structure, with set $\mathcal{H}_1$ of generalized factors.  For $H$
in $\mathcal{H}_1$, let the $\Omega \times \Omega$ relationship matrix
for $H$ be $\mathbf{T}_H$, with corresponding canonical component of
variance $\phi_H$ and common replication $k_H$.  Then 
\[
\sum_{H\in\mathcal{H}_1} \xi_H \mathbf{P}_H = \sum_{H\in\mathcal{H}_1} \phi_H \mathbf{T}_H
\quad \mbox{and} \quad
\xi_H = \sum_{F\in\mathcal{H}_1,\ F\geq H} k_{F} \phi_{F}.
\]
Thus, when both structures are poset block structures,
equation~\eqref{eq:twocov} becomes 
\begin{equation}
\mathbf{V} = \sum_{H\in\mathcal{H}_1} \phi_H \mathbf{T}_H  +
\sum_{H\in\mathcal{H}_2} \psi_H \mathbf{S}_H .
\label{eq:twocancomp}
\end{equation}

As noted in Section~\ref{sec:pbs2}, even for poset block structures
the randomization-based model for variance differs from a 
variance-components model. In equation~\eqref{eq:twocov}, it is the
coefficients $\xi_{\mathbf{P}}$ and $\eta_{\mathbf{Q}}$ which must be nonnegative;
the corresponding canonical components may well be negative,
except for $\phi_{\Omega}$ and $\psi_{\Upsilon}$, which must be positive.

\section{Treatment decomposition and structure balance}
\label{sec:treatsb}

\subsection{Families of expectation models in a two-tiered experiment}
\label{sec:fix}

Consider the two-tiered set-up in Section~\ref{sec:rand2}.
The design function $h$  embeds a copy $V_\Gamma^h$ of
$V_\Gamma$ inside $V_\Upsilon$.
Let $\mathbf{D}_h$ be the $\Gamma \times \Gamma$ diagonal matrix
of replications of treatments.
Then
$\mathbf{X}'_h \mathbf{X}_h = \mathbf{D}_h$,
and the matrix of orthogonal projection onto $V_\Gamma^h$ is
$\mathbf{X}_h \mathbf{D}_h^{-1} \mathbf{X}'_h$, which we 
write as $\mathbf{I}_\mathcal{R}$, because we always associate a
structure $\mathcal{R}$ with $\Gamma$. 
The elements of $\mathcal{R}$ are derived from a family
$\mathcal{M}$ of expectation models on $\Gamma$, as we now show.

With treatment effects fixed, data analysis usually proceeds by
selecting a model from $\mathcal{M}$ and
then estimating the parameters of the chosen model:
see \citep{Bailey08}.
We assume that $\mathcal{M}$
defines an orthogonal decomposition of $V_\Gamma$,
in the following sense.  There is a collection $\mathcal{R}$ of
$\Gamma \times \Gamma$ symmetric, mutually orthogonal, idempotent
matrices whose sum is the $\Gamma \times \Gamma$ identity matrix 
$\mathbf{I}_\Gamma$; each nonzero model in $\mathcal{M}$ is the subspace of
$V_\Gamma$ corresponding to a sum of one or more of the idempotents
in $\mathcal{R}$; if $M$ is such a model then there is at least one
idempotent $\mathbf{R}$ in $\mathcal{R}$ such that $\image(\mathbf{R})
\leq M$ and $M \cap (\image(\mathbf{R}))^\perp$
is in $\mathcal{M}$; each $\mathbf{R}$ in $\mathcal{R}$
occurs at least once in this way, so that it corresponds to
the extra sum of squares for fitting a larger model compared to a
smaller model.

For $\mathbf{R}$ in $\mathcal{R}$, the subspace $\image(\mathbf{R})$
of $V_\Gamma$ is  translated by $h$ into a subspace of $V_\Gamma^h$
whose $\Upsilon \times \Upsilon$ matrix $\mathbf{R}^h$ of orthogonal
projection is $\mathbf{X}_h
\mathbf{R} (\mathbf{R} \mathbf{D}_h \mathbf{R})^- \mathbf{R}
\mathbf{X}_h'$.  We require that $h$ have the property that all such
matrices commute with each other.
When $\mathcal{M}$ is defined by a collection of orthogonal factors on
$\Gamma$, this requirement is equivalent to the condition that the
factors remain orthogonal when considered as factors on $\Upsilon$. 
In the two-tiered context, we shall write $\mathbf{R}^h$ and
$\{\mathbf{R}^h: \mathbf{R}\in\mathcal{R}\}$ simply 
as $\mathbf{R}$ and $\mathcal{R}$ from now on,
so that $\sum_{\mathbf{R}\in\mathcal{R}} \mathbf{R} = \mathbf{I}_\mathcal{R}$.

There is no requirement for the design $h$ to be equireplicate.
For example, suppose that $\Gamma$ consists of the two levels of
a treatment factor. 
If we parameterize the expectations as $\mu+\alpha$ and
$\mu-\alpha$ then the estimators are not orthogonal unless the levels
are equally replicated.  Choice of parametrization should not affect
model fitting, so we prefer to have one model $M_1$ in which
we parameterize the expectations as $\alpha_1$ and $\alpha_2$, with a
submodel $M_2$ in which they are both $\mu$, and a further submodel
$M_3$ in which both expectations are zero.   Then we do have
orthogonality, with $\mathcal{R} = \setof{\mathbf{R}_1, \mathbf{R}_0}$,
where $\mathbf{R}_0
= \card{\Upsilon}^{-1}\mathbf{J}$, which is the projector for the grand
mean, and $\mathbf{R}_1 = \mathbf{I}_\mathcal{R} - \mathbf{R}_0$.

\subsection{Structure balance in a two-tiered experiment}
\label{sec:sb2}

Until Section~\ref{sec:comk} inclusive, 
we insist that $h$ be such that $\mathcal{R}$ is
structure balanced in relation to $\mathcal{Q}$, in the sense defined
in \citep{Brien092}.
This means that there are scalars $\lambda_{\mathbf{QR}}$, for $\mathbf{Q}$
in~$\mathcal{Q}$ and $\mathbf{R}$ in~$\mathcal{R}$, such that
$\mathbf{R}\mathbf{Q}\mathbf{I}_\mathcal{R} = \lambda_{\mathbf{QR}}
\mathbf{R}$.   This equation means that
(i)~$\mathbf{R}\mathbf{Q}\mathbf{R} = \lambda_{\mathbf{QR}} \mathbf{R}$
and 
(ii)~if $\mathbf{R}_1 \ne \mathbf{R}_2$
then $\mathbf{R}_1\mathbf{Q}\mathbf{R}_2 = \mathbf{0}$.
The scalars $\lambda_{\mathbf{QR}}$ are called \emph{efficiency factors}.
It follows that each $\mathbf{Q}$ is the sum of the following mutually
orthogonal idempotents:
(i) $\mathbf{Q} \combine \mathbf{R}$, for all $\mathbf{R}$
in $\mathcal{R}$ with $\lambda_{\mathbf{QR}} \ne 0$, and (ii) if it is
nonzero, $\mathbf{Q} \resid \mathcal{R}$.
These idempotents are
defined by $\mathbf{Q} \combine \mathbf{R} = \lambda_{\mathbf{QR}}^{-1}
\mathbf{Q}\mathbf{R}\mathbf{Q}$ and $\mathbf{Q} \resid \mathcal{R} = \mathbf{Q}
- \sum'_{\mathbf{R}\in\mathcal{R}} \mathbf{Q} \combine \mathbf{R}$, where the
summation is over those $\mathbf{R}$ for which $\lambda_{\mathbf{QR}} \ne 0$.
This set of idempotents is denoted
$\mathcal{Q} \combine \mathcal{R}$
in \citep{Brien092}.

For each $\mathbf{R}$ in $\mathcal{R}$, the efficiency factors
$\lambda_{\mathbf{QR}}$ are nonnegative and sum to~$1$.
If each $\mathbf{R}$ has
some $\mathbf{Q}$ in $\mathcal{Q}$ such that $\lambda_{\mathbf{QR}}=1$
then the structure $\mathcal{R}$ is said to be \emph{orthogonal} 
in relation to the structure $\mathcal{Q}$.

Note that the matrices~$\mathbf{Q}$ are determined by the group of
permutations used for randomizing, and hence ultimately by the
relevant information about $\Upsilon$, such as blocks or managerial
constraints. On the other hand, the matrices~$\mathbf{R}$
depend on the family of expectation models chosen as appropriate.  The
former cannot be altered, but the latter may be refined,
perhaps using pseudofactors, in order to achieve structure balance
\citep{Monod92,Yates36}.
This is achieved by judicious replacement of some matrices
$\mathbf{R}$ in $\mathcal{R}$ by sub-idempotents so that there is
a refinement of the decomposition given by $\mathcal{R}$ into smaller
subspaces:  see  \citep{Brien092}, Section~4. 
Thus we have a larger collection $\mathcal{R}^*$ of mutually
orthogonal idempotents, 
such that each $\mathbf{R}$ in $\mathcal{R}$ is a sum of one or more of the
idempotents in $\mathcal{R}^*$. For example, in a balanced lattice square
design for $k^2$ treatments in $(k+1)/2$ squares,
where $k$~is odd, $\mathcal{R} =
\setof{\mathbf{R}_0, \mathbf{R}_\mathrm{T}}$ where $\mathbf{R}_0$ and
$\mathbf{R}_\mathrm{T}$ are the idempotents for the Mean and
Treatments, respectively.  However, $\mathcal{R}$ is not structure
balanced in relation to the structure $\mathcal{Q}$ defined by
$(k+1)/2$ squares, 
each formed by $k$~rows crossed with $k$~columns.  We form $\mathcal{R}^*$
by replacing $\mathbf{R}_\mathrm{T}$ by $\mathbf{R}_\mathrm{T,\,R}$ and
$\mathbf{R}_\mathrm{T,\,C}$, where these are the idempotents corresponding
to the treatment subspaces partly confounded with rows and columns,
respectively: then $\mathcal{R}^*$ is structure balanced in relation
to $\mathcal{Q}$ (see \citep{Brien092}, Example 5). 

\subsection{Treatment structure and structure balance in a
  three-tiered experiment with two randomizations in a chain} 

Now consider the three-tiered set-up in Section~\ref{sec:derive}.
As in Section~\ref{sec:fix}, the effects $\boldsymbol{\tau}$ are taken
to be fixed, and so we assume that the family of expectation models
gives a set of mutually orthogonal idempotents $\mathbf{R}$
in $\mathcal{R}$ whose sum is the orthogonal projector onto
$V_\Gamma^h$ in $V_\Upsilon$.
Let $\mathbf{M}^h$ be the $\Upsilon \times \Upsilon$ idempotent for
one of these expectation models.  The corresponding $\Omega \times
\Omega$ idempotent $\left(\mathbf{M}^h\right)^f$ is given by
$\left(\mathbf{M}^h\right)^f = \mathbf{X}_f \mathbf{M}^h \left(
    \mathbf{M}^h  \mathbf{D}_f \mathbf{M}^h \right)^- \mathbf{M}^h
    \mathbf{X}_f'
= r^{-1}\mathbf{X}_f
   \mathbf{M}^h \mathbf{X}_f'$
since $\mathbf{D}_f = r^{-1} \mathbf{I}_\Upsilon$.
Therefore, putting $\mathbf{R}^f =
r^{-1}\mathbf{X}_f \mathbf{R} \mathbf{X}_f'$ for $\mathbf{R}$
in~$\mathcal{R}$, we see that the mutually orthogonal idempotents
$\mathbf{R}$ in $\mathcal{R}$ translate to mutually orthogonal idempotents
$\mathbf{R}^f$ on $V_\Omega$.  That is, because $f$~is equireplicate,
the same formula is used to convert both the expectation idempotents
and the variance idempotents from $\Upsilon \times \Upsilon$
matrices to $\Omega \times \Omega$ matrices.
There is still no need for $h$ to be equireplicate.
As in \citep{Brien092}, we shall write $\mathbf{R}^f$ as $\mathbf{R}$
and $\setofall{\mathbf{R}^f}{\mathbf{R}\in \mathcal{R}}$ as $\mathcal{R}$
in the three-tiered context.
We continue to write $\sum_{\mathbf{R} \in \mathcal{R}} \mathbf{R}$ as
$\mathbf{I}_\mathcal{R}$, which is now an $\Omega \times \Omega$ matrix.

In addition to the condition that $f$~be equireplicate, we assume
until Section~\ref{sec:comk} inclusive that
\begin{enumeroman}
\item $\mathcal{Q}$ is
structure balanced in relation to $\mathcal{P}$,
or can be made so, in the sense explained
in Section~\ref{sec:idem};
\item $\mathcal{R}$ is structure balanced 
in relation to $\mathcal{Q}$.
\end{enumeroman}
Then $\mathcal{R}$ is structure balanced
in relation to $\mathcal{P} \combine \mathcal{Q}$, 
$\mathcal{Q} \combine \mathcal{R}$ is
structure balanced in relation to $\mathcal{P}$, and 
$(\mathcal{P} \combine \mathcal{Q}) \combine \mathcal{R} =
\mathcal{P} \combine (\mathcal{Q} \combine \mathcal{R})$,
as shown in \citep{Brien092}.

Let $\mathcal{Q}_1$ be the set of $\mathbf{Q}$ in $\mathcal{Q}$ for which
there is an idempotent $\mathbf{P}$ in $\mathcal{P}$ with
$\lambda_{\mathbf{PQ}}=1$. Define the function~$c$ from
$\mathcal{Q}_1$ to $\mathcal{P}$ such that  $c(\mathbf{Q}) =
\mathbf{P}$ for $\lambda_{\mathbf{PQ}}=1$. 
If $\mathbf{Q} \in \mathcal{Q}_1$ then $\image(\mathbf{Q})  \leq
\image (c(\mathbf{Q}))$,  and
$\mathbf{Q}\mathbf{P} =\mathbf{P}\mathbf{Q} = \mathbf{Q} =
\mathbf{P} \combine \mathbf{Q}$ if $\mathbf{P} = c(\mathbf{Q})$,
while $\mathbf{Q}\mathbf{P}  = \mathbf{P}\mathbf{Q} = \mathbf{0}$
otherwise.  Thus $\mathcal{Q}$ is orthogonal in relation to 
$\mathcal{P}$ when $\mathcal{Q}_1=\mathcal{Q}$.

For $\mathbf{P}\in\mathcal{P}$, equation~\eqref{eq:twocov} shows that
$\mathbf{V}(\mathbf{P}\resid\mathcal{Q}) =
 \xi_{\mathbf{P}}(\mathbf{P}\resid\mathcal{Q})$,
because $(\mathbf{P}\resid\mathcal{Q}) \mathbf{Q} = \mathbf{0}$ for all
$\mathbf{Q}$ in $\mathcal{Q}$.  Hence 
$\image(\mathbf{P}\resid\mathcal{Q})$ is contained in 
an eigenspace of $\mathbf{V}$ with eigenvalue $\xi_{\mathbf{P}}$.
Moreover, if $\mathbf{Q} \in \mathcal{Q}$ and $\lambda_{\mathbf{PQ}}\ne 0$, then
\[
\mathbf{V} (\mathbf{P} \combine \mathbf{Q})
=
\mathbf{V} \frac{\mathbf{PQP}}{\lambda_{\mathbf{PQ}}} =
\frac{\xi_{\mathbf{P}}}{\lambda_{\mathbf{PQ}}} \mathbf{PQP}+
r \sum_{\mathbf{Q}^*} \frac{\eta_{\mathbf{Q}^*}}{\lambda_{\mathbf{PQ}}}
\mathbf{Q}^* \mathbf{PQP}=
\xi_{\mathbf{P}} (\mathbf{P} \combine \mathbf{Q}) + r\eta_{\mathbf{Q}} \mathbf{QP}.
\]
If $\mathbf{Q}\in\mathcal{Q}_1$ and $\mathbf{P} = c(\mathbf{Q})$
then $\image(\mathbf{P} \combine \mathbf{Q})$ is contained in 
an eigenspace of $\mathbf{V}$ with eigenvalue $\xi_{\mathbf{P}} +
r\eta_{\mathbf{Q}}$.  Otherwise, $\image(\mathbf{P} \combine
\mathbf{Q})$ is not contained in any eigenspace of~$\mathbf{V}$. 

If $\mathcal{Q}$ is orthogonal in relation to $\mathcal{P}$ then
\begin{equation}
\mathbf{V} = \sum_{\mathbf{Q}\in\mathcal{Q}}(\xi_{c(\mathbf{Q})} +
r\eta_{\mathbf{Q}})\mathbf{Q} +
\sum_{\mathbf{P}\in\mathcal{P}}\xi_{\mathbf{P}}(\mathbf{P} \resid \mathcal{Q}).
\label{eq:OVS}
\end{equation}
The idempotents in this expression are those in
$\mathcal{P} \combine \mathcal{Q}$, and the the image of each is
contained in an eigenspace of $\mathbf{V}$.  Thus
the set of all positive semidefinite (p.s.d) matrices of the 
form~\eqref{eq:OVS} commute with each other, and have common
eigenspaces: we call this \emph{commutative variance structure} (CVS).
If, in addition,  there is no linear dependence among the coefficients
in \eqref{eq:OVS}, we have OVS.

\subsection{Choice of idempotents}
\label{sec:idem}

The matrices~$\mathbf{P}$ are defined by the group~$G_1$
of permutations used to randomize the design~$f$. The
matrices~$\mathbf{Q}$ are first defined as matrices on $V_\Upsilon$
by the group~$G_2$ of permutations used to randomize the design~$h$,
and then translated by $f$ to matrices on $V_\Omega$. The
matrices~$\mathbf{R}$ depend initially on the chosen family of
expectation models, and are  translated by $h$ and then by~$f$.

Strictly speaking, there is no freedom of choice over the
$\mathbf{Q}$ matrices.  However,
as already outlined for design~$h$ in Section~\ref{sec:sb2}, it is
sometimes possible to turn a design~$f$ without structure balance into
one with structure balance by judicious replacement of some matrices
$\mathbf{Q}$ in $\mathcal{Q}$ by sub-idempotents,
yielding $\mathcal{Q}^*$. The variance matrix in 
equation~\eqref{eq:twocov} is defined by the original $\mathcal{Q}$:
when it is rewritten in terms of $\mathcal{Q}$* it has the constraint
that if $\mathbf{Q}$ in $\mathcal{Q}$ is the sum $\mathbf{Q}_1^* +
\cdots + \mathbf{Q}_n^*$ with $\mathbf{Q}_i^*$ in $\mathcal{Q}^*$ then
each of $\mathbf{Q}_1^*$, \ldots, $\mathbf{Q}_n^*$ has the \emph{same}
spectral component $\eta_\mathbf{Q}$.

There are two types of multiple randomization that form a chain 
as shown in Figure~\ref{fig:samedir}:
see \citep{BrBa:mult,Brien092}. For
\emph{composed} randomizations, the randomizations may be done in
either order, because neither needs knowledge of the outcome of the other.
In contrast, \emph{randomized-inclusive} randomizations have the 
complication that knowledge of the outcome of the randomization of
$\Gamma$ to $\Upsilon$ is needed before $\Upsilon$ can be randomized
to $\Omega$.

As explained in \citep{BrBa:mult}, Section 5.1  and \citep{Brien092}, Section 6,
this knowledge is needed in the second case
partly because the structure~$\mathcal{Q}$ on~$\Upsilon$ defined by the
randomization of design~$h$ is not structure balanced
in relation to~$\mathcal{P}$. Thus $\mathcal{Q}$ needs
to be refined into sub-idempotents or pseudosources, as described above.
The second necessary ingredient for randomized-inclusive randomizations
is that there is at least one source on $\Gamma$
that is confounded or partly confounded with one of the sources on~$\Upsilon$
that needs to be split up. In order to work out the partial confounding of 
sources on $\Gamma$ with those on $\Omega$, it is necessary to keep track of
the partial confounding of the former with the pseudosources on $\Upsilon$.
This may require pseudosources on $\Gamma$.  Most importantly, the unrandomized
version of $f$ is constrained to ensure the correct partial confounding of
(pseudo)sources on $\Gamma$ with those on~$\Omega$.

Although this makes the
procedure more complicated than that for composed randomizations,
the randomization-based model is virtually the same.  As above,
we have to keep track of pseudosources. For
the pseudosources on $\Upsilon$, it is important to remember that pseudosources
of the same source have the \emph{same} spectral component $\eta$.  This
complication can occur for experiments with two composed randomizations when
the second randomization is not consonant: see Example~\ref{eg:cotton}.  It
always occurs for experiments with two randomized-inclusive randomizations:
see Example~\ref{eg:Wheat2}.


\section{Analysis of variance}
\label{sec:sameaov}

\subsection{A two-tiered experiment}
\label{sec:aov}

Consider the two-tiered experiment in Sections~\ref{sec:rand2}
and~\ref{sec:sb2}. In \citep{Brien092} decomposition tables were used
to display the decomposition of $V_{\Upsilon}$
appropriate for such an experiment.  Such a table is a precursor to an
anova table and consists of rows and columns. There is a set of
columns for each tier: one column containing sources, one column
containing degrees of freedom, and, if the design is structure
balanced but not orthogonal, a further column showing efficiency factors.
The sources and pseudosources correspond to idempotents in
$\mathcal{Q}$ or $\mathcal{R}$ which, when they are based on
generalized (pseudo)factors, are labelled as described in Section~3 of
\citep{Brien092}. Each row of the decomposition table corresponds to a
subspace in the decomposition specified by $\mathcal{Q} \combine
\mathcal{R}$. In this paper, we add a column for expected mean squares
to decomposition tables in order to form skeleton anova tables.

The anova table for the analysis of a response variable when the 
variance matrix has the form~\eqref{eq:cov},
the $\tau_i$ are fixed effects and $\mathcal{R}$ is structure balanced
in relation to $\mathcal{Q}$
is given in \citep{Bailey81,Houtman83,Nelder65b}.
The data vector~$\mathbf{y}$ is projected onto each stratum in turn
and then $\mathbf{Qy}$, which is
the projection into stratum $\image(\mathbf{Q})$, is further
decomposed according to the elements 
of $\mathcal{Q} \combine \mathcal{R}$ involving $\mathbf{Q}$.
The following hold.
\begin{enumerate}
\item[(A.a)]
The projections onto different strata are uncorrelated.
\item[(A.b)]
Any orthonormal basis for $\image(\mathbf{Q})$ gives uncorrelated
random variables all with variance $\eta_{\mathbf{Q}}$.
\item[(A.c)]
If $\lambda_{\mathbf{QR}}\ne 0$, then the expected mean square
for $\mathbf{Q} \combine \mathbf{R}$ is equal to
    \[\eta_{\mathbf{Q}} +
    \frac{\lambda_{\mathbf{QR}} \boldsymbol{\tau}' \mathbf{X}'_h \mathbf{R}
    \mathbf{X}_h \boldsymbol{\tau}}{\rank(\mathbf{R})}.\]
\item[(A.d)]
If $\mathbf{Q} \resid \mathcal{R}$ is nonzero, then
the expected mean square for $\mathbf{Q} \resid \mathcal{R}$ is equal to
    $\eta_{\mathbf{Q}}$.
\end{enumerate}
For poset block structures, the spectral component $\eta_{\mathbf{Q}}$ can be
expanded using equation~\eqref{eq:zetaphi} to provide expressions for
the expected mean squares in terms of the canonical components.

The expression $\boldsymbol{\tau}'\mathbf{X}_h'\mathbf{R}\mathbf{X}_h
\boldsymbol{\tau}/\rank(\mathbf{R})$ in~(A.c) is a p.s.d.\ 
quadratic form in the parameters $\tau_i$.  If $\mathcal{R}$ is defined
by a poset block structure  on $\Gamma$ then
$\mathbf{R} =\mathbf{R}_F$ for a generalized factor $F$ on $\Gamma$,
just as  $\mathbf{Q} =\mathbf{Q}_H$ in equation~\eqref{eq:covfact}.
In anova tables, this expression is written as $q(F)$.
In particular,
$q_0 = \boldsymbol{\tau}'\mathbf{X}_h' \mathbf{R}_0
\mathbf{X}_h \boldsymbol{\tau}$, where $\mathbf{R}_0 =
\card{\Upsilon}^{-1}\mathbf{J}$.

\subsection{An experiment with two randomizations in a chain}
\label{sec:aov3}
First consider expectations.
If $\mathbf{P}\in\mathcal{P}$ then
$(\mathbf{P} \resid \mathcal{Q}) \mathbf{I}_{\mathcal{Q}} = \mathbf{0}$.
If, further, $\mathbf{Q}\in\mathcal{Q}$ and $\lambda_{\mathbf{PQ}}\ne0$, then
$((\mathbf{P} \combine \mathbf{Q}) \resid \mathcal{R})
\mathbf{I}_{\mathcal{Q}} \mathbf{I}_{\mathcal{R}}= \mathbf{0}$.
Since
$\E{\mathbf{Y}} = \mathbf{X}_f\mathbf{X}_h\boldsymbol{\tau} =
\mathbf{I}_\mathcal{Q} \mathbf{I}_\mathcal{R}
\mathbf{X}_f\mathbf{X}_h\boldsymbol{\tau}$,
it follows that
$\E{(\mathbf{P} \resid \mathcal{Q})\mathbf{Y}}=
\E{((\mathbf{P} \combine \mathbf{Q}) \resid \mathcal{R})\mathbf{Y}}
=\mathbf{0}$.
If, moreover, $\mathbf{R}\in\mathcal{R}$ and $\lambda_{\mathbf{QR}} \ne 0$,
Section 5 of \citep{Brien092} shows that $(\mathbf{P} \combine
\mathbf{Q}) \combine \mathbf{R} = \lambda^{-1}_{\mathbf{PQ}}
\lambda^{-1}_{\mathbf{QR}}\mathbf{PQRQP}$.
Therefore
\[
\E{((\mathbf{P} \combine \mathbf{Q})\combine \mathbf{R})\mathbf{Y}} =
\frac{1}{\lambda_{\mathbf{PQ}}\lambda_{\mathbf{QR}}} \mathbf{PQRQP}
\mathbf{I}_\mathcal{Q} \mathbf{I}_\mathcal{R} \mathbf{X}_f
\mathbf{X}_h \boldsymbol{\tau} =
\mathbf{PQR} \mathbf{X}_f
\mathbf{X}_h \boldsymbol{\tau}.
\]
Hence
\begin{eqnarray*}
(\E{((\mathbf{P} \combine \mathbf{Q})\combine \mathbf{R})\mathbf{Y}})'
\E{((\mathbf{P} \combine \mathbf{Q})\combine \mathbf{R})\mathbf{Y}}
&=&
\boldsymbol{\tau}'\mathbf{X}_h'\mathbf{X}_f' \mathbf{RQPPQR}
\mathbf{X}_f \mathbf{X}_h \boldsymbol{\tau}\\
& = & \lambda_{\mathbf{PQ}}\lambda_{\mathbf{QR}}
\boldsymbol{\tau}'\mathbf{X}_h'\mathbf{X}_f' \mathbf{R}
\mathbf{X}_f \mathbf{X}_h \boldsymbol{\tau}.
\end{eqnarray*}

Consider a fixed $\mathbf{P}$ in $\mathcal{P}$.
Equation~\eqref{eq:twocov} shows that
\begin{eqnarray}
\Cov{\mathbf{PY}} &=& \mathbf{PVP} =
\xi_\mathbf{P} \mathbf{P} + r\sum_{\mathbf{Q}\in\mathcal{Q}} \eta_\mathbf{Q}
\mathbf{PQP} = \xi_\mathbf{P} \mathbf{P} + r\summ_{\mathbf{Q}\in\mathcal{Q}}
\eta_\mathbf{Q} \lambda_{\mathbf{PQ}} \mathbf{P} \combine \mathbf{Q}
\nonumber
\\
& = & \summ_{\mathbf{Q} \in\mathcal{Q}}
(\xi_\mathbf{P} + r\lambda_{\mathbf{PQ}} \eta_\mathbf{Q})
(\mathbf{P} \combine \mathbf{Q}) + \xi_{\mathbf{P}}(\mathbf{P} \resid
\mathcal{Q}).
\label{eq:pdata}
\end{eqnarray}
Here $\summ_{\mathbf{Q} \in \mathcal{Q}}$ denotes summation over
$\mathbf{Q} \in \mathcal{Q}$ with $\lambda_{\mathbf{PQ}} \neq 0$.
The matrices in  equation~\eqref{eq:pdata} are mutually orthogonal idempotents
which sum to $\mathbf{P}$ and have linearly independent coefficients.
Hence they are the projectors onto the eigenspaces of
$\Cov{\mathbf{PY}}$ with nonzero eigenvalues.  Therefore the results
for $\mathbf{Y}$ in Section~\ref{sec:aov} carry over to 
$\mathbf{PY}$ as follows.
\begin{enumerate}
\item[(A.e)]
The projections onto any two different subspaces of the form
$\image(\mathbf{P}\combine\mathbf{Q})$ or
$\image(\mathbf{P}\resid\mathcal{Q})$
are uncorrelated.
\item[(A.f)]
If $\lambda_{\mathbf{PQ}}\ne 0$,
any orthonormal basis for $\image(\mathbf{P}\combine\mathbf{Q})$ gives
uncorrelated random variables all with variance
$\xi_{\mathbf{P}} + r\lambda_{\mathbf{PQ}}\eta_{\mathbf{Q}}$.
\item[(A.g)]
Any orthonormal basis for $\image(\mathbf{P} \resid \mathcal{Q})$
gives uncorrelated random variables all with variance $\xi_{\mathbf{P}}$.
\item[(A.h)]
If $\lambda_{\mathbf{PQ}}\lambda_{\mathbf{QR}} \ne0$, then
the expected mean square for
$(\mathbf{P} \combine \mathbf{Q}) \combine \mathbf{R}$ is
\[
\xi_\mathbf{P} + r\lambda_{\mathbf{PQ}}\eta_\mathbf{Q} +
\frac{\lambda_{\mathbf{PQ}}\lambda_{\mathbf{QR}}
\boldsymbol{\tau}'\mathbf{X}_h'\mathbf{X}_f'
\mathbf{R} \mathbf{X}_f\mathbf{X}_h\boldsymbol{\tau}}{\rank(\mathbf{R})}.
\]
\item[(A.i)]
If $\lambda_{\mathbf{PQ}}\ne 0$
and $(\mathbf{P} \combine \mathbf{Q}) \resid \mathcal{R}$
is nonzero, then the expected mean square for
$(\mathbf{P} \combine \mathbf{Q}) \resid \mathcal{R}$
is $\xi_\mathbf{P} + r\lambda_{\mathbf{PQ}} \eta_\mathbf{Q}$.
\item[(A.j)]
If $\mathbf{P} \resid \mathcal{Q}$ is nonzero, then
the expected mean square for
$\mathbf{P} \resid \mathcal{Q}$
is   $\xi_\mathbf{P}$.
\end{enumerate}
For poset block structures, the spectral components $\xi_\mathbf{P} $ and $\eta_{\mathbf{Q}}$
can be expanded to express the expected mean squares
in terms of the canonical components.

We write the expression $\boldsymbol{\tau}'\mathbf{X}_h'\mathbf{X}_f'
\mathbf{R}\mathbf{X}_f\mathbf{X}_h\boldsymbol{\tau}/\rank(\mathbf{R})$
as $q(F)$ if $\mathcal{R}$ is defined by a poset block structure on
$\Gamma$ and $\mathbf{R} =\mathbf{R}_F$ for some
generalized factor $F$ on  $\Gamma$.

Finally, consider the whole of $\mathcal{P}$.
If $\mathcal{Q}$ is orthogonal in relation to  $\mathcal{P}$ then we have CVS
and so the projected data corresponding to any two different rows of the
anova table are uncorrelated.   Otherwise, we have the situation, such
as the one in \cite{Ojima98},  where some subspaces corresponding to
idempotents of the form $\mathbf{P}\combine\mathbf{Q}$ do not consist
of eigenvectors of $\mathbf{V}$. Then it is still possible to do anova
in the sense of decomposing the sum of squares of the responses
according to the subspaces, and equating the observed values of the
mean squares to their expectations, but this may not have all the
properties of classical anova.

In particular, let $\mathbf{Q}$ be an idempotent in $\mathcal{Q}$ for
which there are distinct $\mathbf{P}_1$ and $\mathbf{P}_2$ in $\mathcal{P}$
with $\lambda_{\mathbf{P}_1\mathbf{Q}}$ and $\lambda_{\mathbf{P}_2\mathbf{Q}}$ both nonzero. 
Then the projections of the data onto  $\image{(\mathbf{P}_1)}$ and
$\image(\mathbf{P}_2)$ are not independent, because
\begin{eqnarray*}
\Cov{ (\mathbf{P}_1 \combine \mathbf{Q}) \mathbf{Y},
(\mathbf{P}_2 \combine \mathbf{Q}) \mathbf{Y} } &=&
(\mathbf{P}_1 \combine \mathbf{Q})'
\left(
\sum_{\mathbf{P}\in\mathcal{P}} \xi_{\mathbf{P}} \mathbf{P} +
r \sum_{\mathbf{Q}^*\in\mathcal{Q}} \eta_{\mathbf{Q}^*}\mathbf{Q}^*
\right)
(\mathbf{P}_2 \combine \mathbf{Q})\\
&=&
\frac{\mathbf{P}_1 \mathbf{Q} \mathbf{P}_1}{\lambda_{\mathbf{P}_1\mathbf{Q}}}
\left(
\sum_{\mathbf{P}\in\mathcal{P}} \xi_{\mathbf{P}} \mathbf{P} +
r \sum_{\mathbf{Q}^*\in\mathcal{Q}} \eta_{\mathbf{Q}^*}\mathbf{Q}^*
\right)
\frac{\mathbf{P}_2 \mathbf{Q} \mathbf{P}_2}{\lambda_{\mathbf{P}_2\mathbf{Q}}}
\\ &=&
\frac{r\eta_{\mathbf{Q}}}{\lambda_{\mathbf{P}_1\mathbf{Q}}
\lambda_{\mathbf{P}_2 \mathbf{Q}}}
\mathbf{P}_1 \mathbf{Q} \mathbf{P}_1 \mathbf{Q}
\mathbf{P}_2 \mathbf{Q} \mathbf{P}_2 =
r\eta_{\mathbf{Q}} \mathbf{P}_1 \mathbf{Q} \mathbf{P}_2,
\end{eqnarray*}
which has the same rank as $\mathbf{Q}$.
A similar calculation shows that $(\mathbf{P}_1 \combine
\mathbf{Q}_1)\mathbf{Y}$ is not correlated
with $(\mathbf{P}_2 \combine \mathbf{Q}_2)\mathbf{Y}$ 
if $\mathbf{Q}_1 \ne \mathbf{Q}_2$.


\section{Examples}
\label{sec:eg}

Our first two examples show how straightforward the anova 
table is when both designs are orthogonal,
in the sense defined in Section~\ref{sec:sb2}.
Sub\-sequent examples illustrate the application of our results to
other structure-balanced experiments.
Further examples are available in \citep{Brien11}.

\begin{egg}[Meat loaves]
\label{eg:tbb}
The two-phase sensory experiment in Figure~\ref{fig:tbb} is from
\citep{BrBa:mult}, Section 4.1; the design in the second phase
consists of a pair of $6\times 6$ Latin squares in each session.

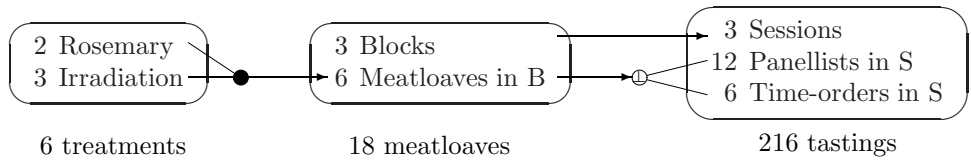
\begin{figure}[htbp]
\setlength{\unitlength}{1cm}
\small
\centering
\begin{picture}(13.3,2)(-0.5,-0.6)
\put(-0.3,0.5){\begin{tierbox}2&Rosemary\\3&Irradiation\end{tierbox}}
\put(3.7,0.5){\begin{tierbox}3&Blocks\\6&Meatloaves in B\end{tierbox}}
\put(8.7,0.5){\begin{tierbox}3&Sessions\\12&Panellists in S\\
6&Time-orders in S\end{tierbox}}
\put(2.9,0.4){\blob}
\put(2.9,0.4){\vector(1,0){1.15}}
\put(2.9,0.4){\line(-1,0){0.7}}
\put(2.9,0.4){\line(-3,2){0.7}}
\put(7.1,0.4){\vector(1,0){1}}
\put(8.2,0.4){\orthcircle}
\put(8.3,0.43){\line(4,1){0.8}}
\put(8.3,0.37){\line(4,-1){0.8}}
\put(7.1,0.95){\vector(1,0){2}}
\put(10.7,-0.5){\makebox(0,0){216 tastings}}
\put(5.4,-0.5){\makebox(0,0){18 meatloaves}}
\put(1.2,-0.5){\makebox(0,0){6 treatments}}
\end{picture}
\caption{Randomization diagram for Example~\ref{eg:tbb}: treatments
  are randomized to meatloaves, which are in turn randomized to
  tastings;
$\mathrm{B}$ denotes Blocks, $\mathrm{S}$ denotes Sessions.}
\label{fig:tbb}
\end{figure}

Table~\ref{tab:tbb}
expands Table~2 of \cite{Brien092} to give
the skeleton anova that
includes the expected mean squares under randomization; there is no
need to show efficiency factors because both designs are orthogonal.

One consequence of this simple orthogonality,
and the lack of pseudosources, is that each
$\eta$-coefficient appears in the final column in conjunction with
exactly one $\xi$-coefficient.
Under randomization all of these coefficients must be nonnegative.
However, canonical components such as $\phi_{\mathrm{SP}}$ and
$\phi_{\mathrm{ST}}$ can be negative, allowing for correlations within
panellists ($\phi_{\mathrm{SP}}$), or within 
time-orders ($\phi_{\mathrm{ST}}$), to be negative.

The appropriate `Residual' for each of the three treatments
sources is the one with $10$~degrees of freedom, which is $ \left(
(\mathrm{P}\inter\mathrm{T}\nesting{\mathrm{S}}) \combine
(\mathrm{Meatloaves}\nesting{\mathrm{B}}) \right) \resid \mathcal{R} $,
where $\mathcal{R}$ is the structure on the treatments tier.

\begin{table}
\caption{\label{tab:tbb}Skeleton analysis of variance for
Example~\ref{eg:tbb}} 
\begin{center}
\begin{tabular}{*{3}{lrc}l}
\toprule
\multicolumn{2}{c}{\textbf{tastings tier}} &
& \multicolumn{2}{c}{\textbf{meatloaves tier}} & &
\multicolumn{2}{c}{\textbf{treatments tier}} & & \\
\cmidrule{1-2} \cmidrule{4-5} \cmidrule{7-8}
\textbf{source} & \textbf{d.f.} & &
\textbf{source} & \textbf{d.f.} & &
\textbf{source} & \textbf{d.f.} & &
\multicolumn{1}{c}{\textbf{E.M.S.}} \\
\midrule
Mean & 1 & & Mean & 1 & & Mean & 1 & &  $\xi_0 + 12\eta_0 + q_0$\\
\midrule
    Sessions & 2 & & Blocks & 2 & & & & & $\xi_\mathrm{S} + 12\eta_\mathrm{B}$ \\
\midrule
     $\mathrm{Panellists}\nesting{\mathrm{S}}$ 
&33 & & & & & &  & & $\xi_\mathrm{SP}$ \\
\midrule
     $\mathrm{Time\mbox{-}orders}\nesting{\mathrm{S}}$ 
&15 & & & & &  & & & $\xi_\mathrm{ST}$ \\
\midrule
     $\mathrm{P}\inter\mathrm{T}\nesting{\mathrm{S}}$ &165 & &
$\mathrm{Meatloaves}\nesting{\mathrm{B}}$& 15 & & Rosemary & 1 & &
$\xi_\mathrm{STP} + 12\eta_\mathrm{BM} + q(\mathrm{R})$ \\
\cmidrule{7-10}
 & & & & & & Irradiation & 2 & & $\xi_\mathrm{STP} + 12\eta_\mathrm{BM} +
q(\mathrm{I})$\\
\cmidrule{7-10}
& & & & & & $\mathrm{R} \inter \mathrm{I}$ & 2 & & $\xi_\mathrm{STP} +
12\eta_\mathrm{BM} + q(\mathrm{RI})$ \\
\cmidrule{7-10}
& & & & & & Residual & 10 & & $\xi_\mathrm{STP} + 12\eta_\mathrm{BM}$ \\
\cmidrule{4-10}
& & & Residual & 150 & & & & & $\xi_\mathrm{STP}$ \\
\bottomrule
\end{tabular}\end{center}
\end{table}
\end{egg}

\begin{egg}[Cotton fibres]
\label{eg:cotton} Example 4 of \citep{BrBa:mult} describes a two-phase
experiment that consists of a field phase in which cotton is produced
in a field trial and a testing phase in which the cotton is tested for
strength. In this example, the use of both the spectral and canonical
components will be illustrated. 

The randomization diagram for Plan B is in
Figure~\ref{fig:cox4}. As discussed in \citep{BrBa:mult},
there is a lack of
consonance between the randomization of fibres to tests and the
nesting of the associated factors: nested Fibres is randomized to
nesting Operatives and nesting Block$\meet$Plots is randomized to
nested Tests.  This forces us to
introduce a two-level pseudofactor~$\mathrm{F}_1$ that indexes the two
groups of 15 fibres; it neither is nested in nor nests anything.

The middle panel in Figure~\ref{fig:cox4} indicates that the proper
randomization for the first phase is to randomly permute blocks,
randomly permute plots independently within each block, and then
randomly permute fibres within each plot.
As noted in \citep{BrBa:mult}, the two fibres taken from each plot must be
independently randomized to the two levels of the pseudofactor
$\mathrm{F}_1$, so that its levels do not correspond to any
inherent property such as strength or length.

\begin{figure}[htbp]
\setlength{\unitlength}{1cm}
\small
\centering
\begin{picture}(12.6,1.9)(0.4,-0.35)
\put(1.0,0.9){\begin{tierbox}$5$ & K\end{tierbox}}
\put(4,0.85){\begin{tierbox}$2$ & Fibres in B, P \\
                            $5$ & Plots in B\\$3$ & Blocks\end{tierbox}}
\put(9.45,0.85){\begin{tierbox}$2$ & Operatives\\ \\
                              $15$ & Tests in O\end{tierbox}}
\put(2.2,0.95){\vector(1,0){2.15}}
\put(7.9,0.525){\line(-3,1){1}}
\put(7.9,0.525){\blob}
\put(6.9,0.525){\vector(1,0){2.9}}
\put(7.0,1.35){\line(1,0){1.2}}
\put(8.3,1.35){\makebox(1,0)[l]{$2 \hspace{\nlevnamesep} {\rm F}_1$}}
\put(9,1.35){\vector(1,0){0.9}}
\put(1.8,-0.2){\makebox(0,0){5 treatments}}
\put(5.75,-0.2){\makebox(0,0){30 fibres}}
\put(11,-0.2){\makebox(0,0){30 tests}}
\end{picture}
\caption{\label{fig:cox4}Composed randomizations in
  Example~\ref{eg:cotton}: treatments are randomized to fibres, which
are in turn randomized to tests;
K denotes Potash treatments;
$\mathrm{B}$ denotes Blocks, $\mathrm{P}$ denotes Plots;
$\mathrm{O}$ denotes Operatives; $\mathrm{F}_1$ is a pseudofactor for
  Fibres.}
\end{figure}
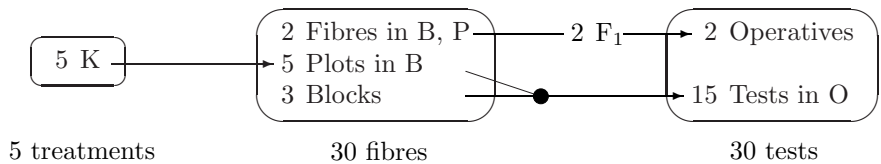

The variance matrix under the randomizations is
\begin{eqnarray*}
\mathbf{V} & = &
    \xi_0 \mathbf{P}_0 + \xi_\mathrm{O} \mathbf{P}_\mathrm{O}
    + \xi_\mathrm{OT} \mathbf{P}_\mathrm{OT} + \eta_0 \mathbf{Q}_0 +
      \eta_\mathrm{B} \mathbf{Q}_\mathrm{B} + \eta_\mathrm{BP}
      \mathbf{Q}_\mathrm{BP} + \eta_\mathrm{BPF} \mathbf{Q}_\mathrm{BPF} \\
    & = &  \phi_{0} \mathbf{T}_0 + \phi_\mathrm{O} \mathbf{T}_\mathrm{O}
    + \phi_\mathrm{OT} \mathbf{T}_\mathrm{OT} + \psi_{0} \mathbf{S}_0 +
      \psi_\mathrm{B} \mathbf{S}_\mathrm{B} + \psi_\mathrm{BP}
      \mathbf{S}_\mathrm{BP} + \psi_\mathrm{BPF} \mathbf{S}_\mathrm{BPF},
\end{eqnarray*}
where $\mathbf{T}_H$ and $\mathbf{S}_H$ are the relationship matrices
and $\phi_{H}$ and $\psi_{H}$ are the canonical components for a
generalized factor $H$ from the tests or fibres tiers, respectively.

There is no $\mathbf{Q}$-matrix
for the pseudofactor: it is irrelevant to the randomization of
treatments to fibres and is not one of the unrandomized factors,
that gives rise to covariance, in the other randomization.
However, we can rewrite $\mathbf{Q}_\mathrm{BPF}$ as the sum of two
$\mathbf{Q}^*$ matrices, one for each of $\mathrm{F}_1$ and
$\mathrm{Fibres}\nesting{\mathrm{Plots}\wedge\mathrm{Blocks}} \resid
\mathrm{F}_1$; the coefficient of both is $\eta_\mathrm{BPF}$.

The following expressions show how the canonical components in this
example measure excess covariance: 
\begin{equation*}
\begin{array}{l}
  \phi_\mathrm{OT} = \gamma_\mathrm{OT} - \gamma_\mathrm{O},\;
      \phi_\mathrm{O} = \gamma_\mathrm{O} - \gamma_{0},\;
      \phi_{0} = \gamma_{0},\; \\
  \psi_\mathrm{BPF} = \zeta_\mathrm{BPF} - \zeta_\mathrm{BP},\;
     \psi_\mathrm{BP} = \zeta_\mathrm{BP} - \zeta_\mathrm{B},\;
     \psi_\mathrm{B} = \zeta_\mathrm{B} - \zeta_{0},\;
     \psi_{0} = \zeta_{0}.\; \\
\end{array}
\end{equation*}
Thus, $\phi_{0}$, $\phi_\mathrm{O}$ and $\phi_\mathrm{OT}$ measure,
respectively, the basic covariance of `unrelated' tests, the excess of
the covariance of different tests by the same operator over that of
`unrelated' tests, and the excess of the 
(co)variance of the same tests over that of different tests by the
same operator. The $\psi$-parameters from the fibres tier can be
similarly interpreted using the $\zeta$-parameters. 

The skeleton anova is in Table~\ref{tab:cox4}; again, there is no
need to show efficiency factors because both designs are orthogonal.
Now the coefficent $\eta_{\mathrm{BPF}}$ occurs with two
  different $\xi$-coefficients.
This is because the source
$\mathrm{Fibres}\nesting{\mathrm{Plots}\wedge\mathrm{Blocks}}$ has
  been split into two by the pseudofactor.

\begin{sidewaystable}
\caption{\label{tab:cox4}Skeleton analysis of variance for
Example~\ref{eg:cotton} with expected mean squares in terms of
spectral components and canonical components} 
\begin{center}
\begin{tabular}{*{3}{lrc}lc*{5}{l@{}}}
\toprule
\multicolumn{2}{c}{\textbf{tests tier}} & &
\multicolumn{2}{c}{\textbf{fibres tier}} & &
\multicolumn{2}{c}{\textbf{treatments tier}} 
& & \multicolumn{7}{c}{\textbf{E.M.S.}} \\
\cmidrule{1-2} \cmidrule{4-5} \cmidrule{7-8}
\textbf{source} & \textbf{d.f.} & &
\textbf{source} & \textbf{d.f.} & &
\textbf{source} & \textbf{d.f.} & &
\multicolumn{1}{c}{\textbf{spectral components}} & &
\multicolumn{5}{c}{\textbf{canonical components}}\\
\midrule
Mean & 1 & & Mean & 1 & & Mean & 1 & & $\xi_0 + \eta_0 + q_0$ & &
         $\phi_\mathrm{OT}$ & $+ 15\phi_\mathrm{O}$ & $+ 30\phi_{0}$ \\
 & & & & & & & & & & & & &
         $+ \psi_\mathrm{BPF}$ & $+ 2\psi_\mathrm{BP}$ & $+ 10\psi_\mathrm{B}$ \\
  & & & & & & & & & & &  & & $+ 30\psi_{0}$ & $+ q_{0}$\\
\midrule
Operatives & 1 & & $\mathrm{F}_1$  & 1 & & & & &
         $\xi_\mathrm{O} + \eta_\mathrm{BPF}$ & &
         $\phi_\mathrm{OT}$ & $+ 15\phi_\mathrm{O}$ & $+ \psi_\mathrm{BPF}$ \\
\midrule
$\mathrm{Tests}\nesting{\mathrm{O}}$& 28 & & Blocks & 2 & & & & &
        $\xi_\mathrm{OT} + \eta_\mathrm{B}$& &
        $\phi_\mathrm{OT}$ & & $+ \psi_\mathrm{BPF}$ & $+
2\psi_\mathrm{BP}$ & $+ 10\psi_\mathrm{B}$\\
\cmidrule{4-16}
& & & $\mathrm{Plots}\nesting{\mathrm{B}}$& 12 & & K & 4 & &
       $\xi_\mathrm{OT} + \eta_\mathrm{BP} + q(\mathrm{K})$& &
       $\phi_\mathrm{OT}$ & & $+ \psi_\mathrm{BPF}$ & $+ 2\psi_\mathrm{BP}$ &
       $ + q(\mathrm{K})$   \\
& & & & & & Residual & 8 & &    $\xi_\mathrm{OT} + \eta_\mathrm{BP}$ & &
       $\phi_\mathrm{OT}$ & & $+ \psi_\mathrm{BPF}$ & $+ 2\psi_\mathrm{BP}$ \\
\cmidrule{4-16}
& & &  $\mathrm{Fibres}\nesting{\mathrm{P}\wedge\mathrm{B}}
\resid \mathrm{F}_1$& 14 & & & & & $\xi_\mathrm{OT} + \eta_\mathrm{BPF}$
   & & $\phi_\mathrm{OT}$ & & $+ \psi_\mathrm{BPF}$\\
\bottomrule
\end{tabular}
\end{center}
\end{sidewaystable}
\end{egg}

\begin{egg}[Two-phase sensory experiment]
\label{eg:2PhaseSensory}
Section~3 of \citep{Brien99} describes a two-phase sensory experiment.
The first, or field, phase is a viticultural experiment and the second, or
evaluation, phase involves the assessment of wine made from the produce
of the first-phase plots.
The randomization diagram for it, given in \citep{BrBa:mult}, is in
Figure~\ref{fig:2PhaseSensory} and the decomposition table for it is
derived in \citep{Brien092}, Example~1. Here the
skeleton anova is in Table~\ref{tab:decompsensory}.
Although there are pseudo\-factors for the Judges factor in the
evaluations tier, they are ignored in doing the randomization as
the six judges are permuted with no distinction. These pseudofactors are
used only to obtain the systematic layout;
they do not give rise to pseudosources.

\begin{figure}[h]
\setlength{\unitlength}{0.9cm}
\footnotesize
\centering
\begin{picture}(13.9,4.0)(-1.6,0.7)
\put(-1.7,2.4){\begin{tierbox}4 & Trellis\vspace{5pt} \\
2 & Methods\end{tierbox}}
\put(-1.3,0.9){\makebox(1,0)[l]{8 treatments}}
\put(0.4,2.8){\vector(1,0){0.7}}
\put(1.2,2.8){\nonorthcircle}
\put(1.3,2.8){\line(4,1){0.4}}
\put(1.3,2.8){\line(6,-1){0.4}}
\put(0.6,2.2){\vector(1,0){1.1}}
\put(1.3,2.7){\begin{tierbox}2 & Squares\\3 & Rows\\
                             4 & Columns in ${\rm Q}$\\
                             2 & Halfplots in ${\rm Q},\
{\rm R},\ {\rm C}$\end{tierbox}}
\put(2.7,0.9){\makebox(1,0)[l]{48 half-plots}}
\put(4.8,3.5){\vector(4,1){3.4}}
\put(5.6,2.1){\vector(4,-1){2.55}}
\put(4.8,3.05){\vector(1,0){1.3}}
\put(6.2,3.05){\orthcircle}
\put(6.32,3.05){\line(5,2){1.75}}
\put(6.32,3.05){\line(1,0){0.3}}
\put(4.8,2.65){\vector(1,0){1.3}}
\put(6.2,2.65){\nonorthcircle}
\put(6.32,2.65){\line(1,0){0.3}}
\put(6.32,2.65){\line(3,-1){1.75}}
\put(6.7,3.05){\makebox(1,0)[l]{$3 \hspace{\nlevnamesep} {\rm J}_2$}}
\put(6.7,2.65){\makebox(1,0)[l]{$2 \hspace{\nlevnamesep} {\rm J}_1$}}
\put(7.65,2.9){\blob}
\put(7.75,2.9){\line(-2,1){0.4}}
\put(7.75,2.9){\line(-2,-1){0.4}}
\put(7.55,2.9){\line(1,0){0.55}}
\put(7.7,2.8){\begin{tierbox} 2 & Occasions\\[5pt]
                              3 & Intervals in ${\rm O}$\\[10pt] 6
& Judges\\[10pt]
                              4 & Sittings in ${\rm O},\ {\rm I}$\\[5pt]
                              4 & Positions in ${\rm O},\ {\rm I},
\ {\rm S},\ {\rm J}$
\end{tierbox}}
\put(8.9,0.8){\makebox(1,0)[l]{576 evaluations}}
\end{picture}
\caption[]{Randomization diagram for Example~\ref{eg:2PhaseSensory}: treatments
  are randomized to half-plots,
  which are, in turn, randomized to
  evaluations; $\mathrm{Q}$, $\mathrm{R}$, $\mathrm{C}$, $\mathrm{O}$,
  $\mathrm{I}$, $\mathrm{S}$, $\mathrm{J}$ denote Squares, Rows,
  Columns, Occasions, Intervals, Sittings and 
  Judges, respectively; ${\rm J}_1$ and ${\rm J}_2$ are pseudofactors
  for Judges.}
\label{fig:2PhaseSensory}
\end{figure}
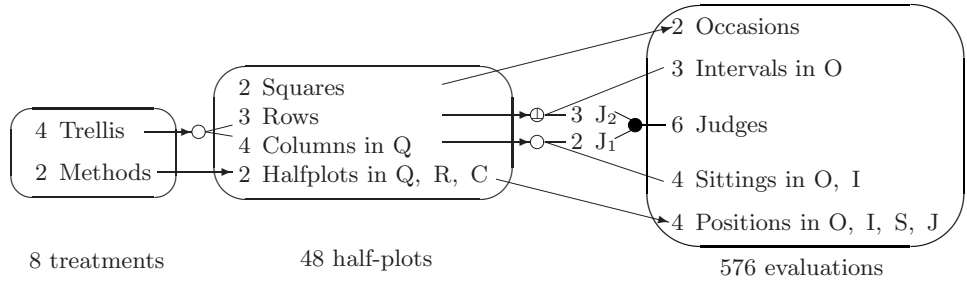

In this example, neither design is orthogonal, and so efficiency
factors need to be shown in the anova table.  The only
source in the halfplots tier which is not orthogonal to the sources
in the evaluations tier is
$\mathrm{Columns}\nesting{\mathrm{Squares}}$: the corresponding
coefficient $\eta_{\mathrm{QC}}$ occurs in conjunction with two
different $\xi$-coefficients.  Similarly, the treatment source Trellis
is nonorthogonal to three idempotents in $\mathcal{P} \combine
\mathcal{Q}$, and so information about Trellis differences is
available in three different subspaces, as shown by the three
occurrences of $q(\mathrm{T})$ in the table.

\begin{table}
\caption{\label{tab:decompsensory}Skeleton analysis of variance table
for Example~\ref{eg:2PhaseSensory} (O = Occasions, I = Intervals, S =
Sittings, J = Judges, P = Positions, Q = Squares, C =  Columns, R =
Rows, H = Halfplots, T = Trellis, M = Methods)}
\begin{center}
\hspace*{-1cm}
\begin{tabular}{lrc*{2}{clrc}l}
\toprule
\multicolumn{2}{c}{\textbf{evaluations tier}} & &
\multicolumn{3}{c}{\textbf{halfplots tier}} & &
\multicolumn{3}{c}{\textbf{treatments tier}} & \\
\cmidrule{1-2} \cmidrule{4-6} \cmidrule{8-10}
\textbf{source} & \textbf{d.f.} & &
\textbf{eff.} & \textbf{source} & \textbf{d.f.} & &
\textbf{eff.} & \textbf{source} & \textbf{d.f.} & &
\multicolumn{1}{c}{\textbf{E.M.S.}} \\
\midrule
Mean & 1 & & 1 & Mean & 1 & & 1 & Mean & 1 & & $\xi_0 + 12\eta_0 + q_0$ \\
\midrule
$\mathrm{O}$ & 1 & & 1 & $\mathrm{Q}$ & 1 & & & & & &
                             $\xi_\mathrm{O} + 12\eta_\mathrm{Q}$ \\
\midrule
$\mathrm{I} \nesting{\mathrm{O}}$ & 4 & & & & & & & & & & $\xi_\mathrm{OI}$ \\
\midrule
$\mathrm{S} \nesting{\mathrm{O} \wedge \mathrm{I}}$ & 18 &
   & \raisebox{8.5pt}[9pt][5pt]{}$\frac{1}{3}$ & $\mathrm{C} \nesting{\mathrm{Q}}$ & 6 &
   & \raisebox{8.5pt}[9pt][5pt]{}$\frac{1}{27}$ & $\mathrm{T}$ & 3 & &
                   $\xi_\mathrm{OIS}
                   +
		   \raisebox{8.5pt}[9pt][5pt]{}\frac{1}{3}12\eta_\mathrm{QC}
+ \raisebox{8.5pt}[9pt][5pt]{}\frac{1}{27}q(\mathrm{T})$ \\
& & & & & & & & Residual & 3 & &
                   $\xi_\mathrm{OIS}
                   + \raisebox{8.5pt}[9pt][5pt]{}\frac{1}{3}12\eta_\mathrm{QC}$ \\
\cmidrule{4-12}
& & & & Residual & 12 & & & & & & $\xi_\mathrm{OIS}$\\
\midrule
\mbox{J} & 5 & & & & & & & & & & $\xi_\mathrm{J}$ \\
\midrule
$\mathrm{O} \inter \mathrm{J}$ & 5 & & & & & & & & & & $\xi_\mathrm{OJ}$ \\
\midrule
$\mathrm{I} \inter \mathrm{J} \nesting{\mathrm{O}}$ & 20 & &
  1 & $\mathrm{R}$ & 2 & & & & & &
                   $\xi_\mathrm{OIJ} + 12\eta_\mathrm{R}$ \\
& & & 1 & $\mathrm{R} \inter \mathrm{Q}$ & 2 & & & & & &
                   $\xi_\mathrm{OIJ} + 12\eta_\mathrm{QR}$ \\
& & & & Residual & 16 & & & & & & $\xi_\mathrm{OIJ}$ \\
 \midrule
$\mathrm{S} \inter \mathrm{J} \nesting{\mathrm{O} \wedge \mathrm{I}}$
& 90 & & \raisebox{8.5pt}[9pt][5pt]{}$\frac{2}{3}$ &
 $\mathrm{C} \nesting{\mathrm{Q}}$ & 6 &
   & \raisebox{8.5pt}[9pt][5pt]{}$\frac{2}{27}$ & T & 3 & &
              $\xi_\mathrm{OISJ}
              + \raisebox{8.5pt}[9pt][5pt]{}\frac{2}{3}12\eta_\mathrm{QC}
              + \raisebox{8.5pt}[9pt][5pt]{}\frac{2}{27}q(\mathrm{T})$ \\
& & & & & & & & Residual & 3 & & $\xi_\mathrm{OISJ}
         + \raisebox{8.5pt}[9pt][5pt]{}\frac{2}{3}12\eta_\mathrm{QC}$ \\
\cmidrule{4-12}
& & & 1 & $\mathrm{R} \inter \mathrm{C} \nesting{\mathrm{Q}}$ & 12 &
  & \raisebox{8.5pt}[9pt][5pt]{}$\frac{8}{9}$ & T & 3 & &
              $\xi_\mathrm{OISJ} + 12\eta_\mathrm{QRC}
              + \raisebox{8.5pt}[9pt][5pt]{}\frac{8}{9}q(\mathrm{T})$ \\
& & & & & & & & Residual & 9 & & $\xi_\mathrm{OISJ} + 12\eta_\mathrm{QRC}$ \\
\cmidrule{4-12}
& & & & Residual & 72 & & & & & & $\xi_\mathrm{OISJ}$\\
\midrule
$\mathrm{P} \nesting{\mathrm{O} \wedge \mathrm{I} \wedge \mathrm{S}
\wedge \mathrm{J}}$ & 432 & & 1 & $\mathrm{H} \nesting{\mathrm{Q} \wedge
\mathrm{R} \wedge \mathrm{C}}$ & 24 &
& 1 & M & 1 & & $\xi_\mathrm{OISJP} + 12\eta_\mathrm{QRCH} + q(\mathrm{M})$ \\
& & & & & & & 1 & $\mathrm{T} \inter \mathrm{M}$ & 3 & &
           $\xi_\mathrm{OISJP} + 12\eta_\mathrm{QRCH} +  q(\mathrm{TM})$  \\
& & & & & & & & Residual & 20 & & $\xi_\mathrm{OISJP} + 12\eta_\mathrm{QRCH}$ \\
\cmidrule{4-12}
& & & & Residual & 408 & & & & & & $\xi_\mathrm{OISJP}$\\
\bottomrule
\end{tabular}
\end{center}
\end{table}
\end{egg}

\begin{egg}[Duplicated wheat measurements]
\label{eg:Wheat2} Example~9 of \citep{BrBa:mult} is an experiment
with a field phase and a laboratory phase. In the field phase 49
lines of wheat are investigated using a randomized complete-block
design with four blocks. Here the laboratory phase is modified by
supposing that the procedure described in \citep{BrBa:mult} is repeated
on a second occasion. That is, two samples are 
obtained from each plot, one to be processed on each occasion.
Figure~\ref{fig:wheat2} gives the randomization
diagram for the modified experiment. Recall that a $7 \times 7$
balanced lattice square design with four replicates is used to assign
the blocks, plots and lines to four intervals in each occasion. In each
interval on each
occasion there are seven runs at which samples are
processed at seven consecutive times. Pseudofactors are introduced for
lines and plots in order to define the design of the second phase.

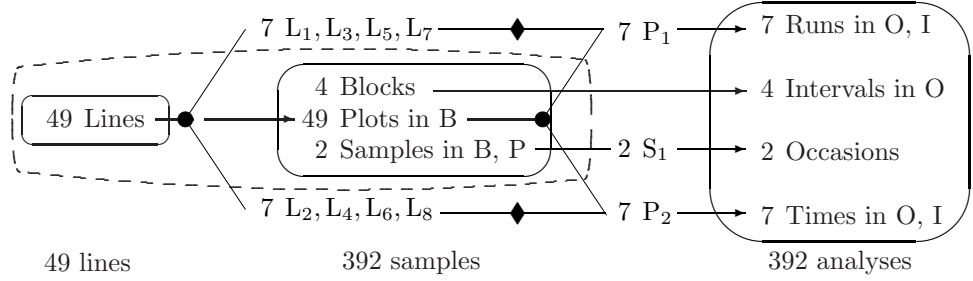
\begin{figure}[htbp]
\setlength{\unitlength}{1cm} \centering \small
\begin{picture}(12.7,3.5)(0,0)
\curvedashes[4pt]{0,1,1} \closecurve(0.1,2.0, 0.2,1.3, 7.7,1.2, 7.8,2.0
,7.7,2.8, 0.2,2.7)
\put(0.1,1.9){\begin{tierbox}49 & Lines\hspace{250pt}\end{tierbox}} \put(2.65,2.0){\vector(1,0){1.2}}
\put(2,2.0){\line(1,0){0.5}}
\put(2.4,2.0){\blob}
\put(2.4,2.0){\line(2,3){0.8}}
\put(2.4,2.0){\line(2,-3){0.8}}
\put(3.4,3.2){\makebox(1,0)[l]{$7 \hspace{\nlevnamesep}
                                {\rm L}_1,{\rm L}_3,{\rm L}_5,{\rm L}_7$}}
\put(3.4,0.8){\makebox(1,0)[l]{$7 \hspace{\nlevnamesep}
                                {\rm L}_2,{\rm L}_4,{\rm L}_6,{\rm L}_8$}}
\put(3.5,1.9){\begin{tierbox}4 & Blocks\\49 & Plots in B\\
                             2 & Samples in B, P\end{tierbox}}
\put(5.7,2.375){\vector(1,0){4.15}}
\put(7.15,2.0){\line(-1,0){1.0}}
\put(7.15,2.0){\blob}
\put(7.15,2.0){\line(2,3){0.805}}
\put(7.15,2.0){\line(2,-3){0.85}}
\put(8.15,3.15){\makebox(1,0)[l]{$7\hspace{\nlevnamesep} {\rm P}_1$}}
\put(8.15,0.75){\makebox(1,0)[l]{$7\hspace{\nlevnamesep} {\rm P}_2$}}
\put(6.7,3.10){\makepseudo}
\put(6.7,0.65){\makepseudo}
\put(7.96,3.2){\line(-1,0){2.2}}
\put(7.99,0.75){\line(-1,0){2.2}}
\put(8.95,3.2){\vector(1,0){0.90}}
\put(8.95,0.75){\vector(1,0){0.90}}
\put(8.15,1.6){\makebox(1,0)[l]{$2\hspace{\nlevnamesep} {\rm S}_1$}}
\put(8.05,1.6){\line(-1,0){1.0}}
\put(8.95,1.6){\vector(1,0){0.90}}
\put(9.25,1.875){\begin{tierbox}\quad 7 & Runs in O, I\\&\\
                                 4 & Intervals in O\\&\\
                                 2 & Occasions\\&\\ 7 & Times in O, I\end{tierbox}}
\put(1.1,0.1){\makebox(0,0){49 lines}}
\put(5.4,0.1){\makebox(0,0){392 samples}}
\put(11.1,0.1){\makebox(0,0){392 analyses}}
\end{picture}
\caption{Randomized-inclusive randomizations in
  Example~\ref{eg:Wheat2}: lines are randomized to
samples, which are in turn
randomized to analyses; $\mathrm{B}$, $\mathrm{P}$, 
$\mathrm{O}$, $\mathrm{I}$ denote Blocks, Plots, 
Occasions, Intervals, respectively;
$\mathrm{L}_1$, \ldots,
  $\mathrm{L}_8$ are mutually orthogonal pseudofactors for Lines;
  $\mathrm{P}_1$ and $\mathrm{P}_2$ are pseudofactors for Plots,
  determined from different Lines pseudofactors in different blocks;
$\mathrm{S}_1$ is a pseudofactor for Samples.}
\label{fig:wheat2}
\end{figure}

\begin{table}
\caption{\label{tab:ANOVAWheat2}Skeleton analysis of variance for
Example~\ref{eg:Wheat2}} 
\begin{center}
\begin{tabular}{*{2}{lrc}clrcl}
\toprule
\multicolumn{2}{c}{\textbf{analyses tier}} & &
\multicolumn{2}{c}{\textbf{samples tier}} &  &
\multicolumn{3}{c}{\textbf{lines tier}} & \\
\cmidrule{1-2} \cmidrule{4-5} \cmidrule{7-9}
\textbf{source} & \textbf{d.f.} & &
\textbf{source} & \textbf{d.f.} & &
\textbf{eff.} & \textbf{source} & \textbf{d.f.} & &
\multicolumn{1}{c}{\textbf{E.M.S.}} \\
\midrule
Mean & 1 & & Mean & 1 & & & Mean & 1 & & $\xi_0 + \eta_0 + q_0$ \\
\midrule
Occasions & 1 & & $\mathrm{S}_1$ & 1 & & & & & &
                          $\xi_\mathrm{O} + \eta_\mathrm{BPS}$ \\
\midrule
$\mathrm{Intervals}\nesting{\mathrm{O}}$ & 6 & & Blocks & 3 & & & & & &
                          $\xi_\mathrm{OI} + \eta_\mathrm{B}$ \\
& & & $\mathrm{S}_1\inter\mathrm{B}$ & 3 & & & & & &
                          $\xi_\mathrm{OI} + \eta_\mathrm{BPS}$ \\
\midrule
$\mathrm{Runs} \nesting{\mathrm{O}\wedge\mathrm{I}}$ & 48 & &
       $\mathrm{P}_1 \nesting{\mathrm{B}}$ &  24 &
        &\raisebox{8.5pt}[9pt][5pt]{}$\frac{1}{4}$ &  $\mathrm{Lines}_{\mathrm{R}}$  & 24 &
        &\raisebox{8.5pt}[9pt][5pt]{}$\xi_\mathrm{OIR} + \eta_\mathrm{BP} + \frac{1}{4}q(L_R)$ \\
& & &  $\mathrm{S}_1\inter\mathrm{P}_1\nesting{\mathrm{B}}$ & 24 & & & & & &
                          $\xi_\mathrm{OIR} + \eta_\mathrm{BPS}$ \\
\midrule
$\mathrm{Times}\nesting{\mathrm{O}\wedge\mathrm{I}}$ & 48 & &
       $\mathrm{P}_2 \nesting{\mathrm{B}}$ &  24 &
        &\raisebox{8.5pt}[9pt][5pt]{}$\frac{1}{4}$ & $\mathrm{Lines}_{\mathrm{T}}$ & 24 &
        &\raisebox{8.5pt}[9pt][5pt]{}$\xi_\mathrm{OIT} + \eta_\mathrm{BP} + \frac{1}{4}q(L_T)$\\
& & &  $\mathrm{S}_1\inter\mathrm{P}_2\nesting{\mathrm{B}}$ & 24 & & & & & &
                           $\xi_\mathrm{OIT} + \eta_\mathrm{BPS}$ \\
\midrule
$\mathrm{R}\inter\mathrm{T}\nesting{\mathrm{O}\wedge\mathrm{I}}$ & 288 & &
       $\mathrm{Plots}\nesting{\mathrm{B}}_{\sresid}$ &  144 &
        &\raisebox{8.5pt}[9pt][5pt]{}$\frac{3}{4}$ & $\mathrm{Lines}_{\mathrm{R}}$ & 24 &
        &\raisebox{8.5pt}[9pt][5pt]{}$\xi_\mathrm{OIRT} + \eta_\mathrm{BP} + \frac{3}{4}q(L_R)$ \\
& & & && & \raisebox{8.5pt}[9pt][5pt]{}$\frac{3}{4}$ &
        $\mathrm{Lines}_{\mathrm{T}}$ & 24 & &
         \raisebox{8.5pt}[9pt][5pt]{}$\xi_\mathrm{OIRT} + \eta_\mathrm{BP} + \frac{3}{4} q(L_T)$\\
& & & & & & & Residual  & 96 & & $\xi_\mathrm{OIRT} + \eta_\mathrm{BP}$  \\
\cmidrule{4-11}
& & &
$\mathrm{Samples}\nesting{\mathrm{B}\wedge\mathrm{P}}_{\sresid}$
              & 144 & & & & & & $\xi_\mathrm{OIRT} + \eta_\mathrm{BPS}$ \\
\bottomrule
\end{tabular}
\end{center}
\end{table}

The variance matrix under the randomizations is
\begin{eqnarray}
\mathbf{V} & = & \xi_0 \mathbf{P}_0 + \xi_\mathrm{O}
\mathbf{P}_\mathrm{O} + \xi_\mathrm{OI} \mathbf{P}_\mathrm{OI} +
\xi_\mathrm{OIR} \mathbf{P}_\mathrm{OIR} + \xi_\mathrm{OIT}
\mathbf{P}_\mathrm{OIT} +
\xi_\mathrm{OIRT} \mathbf{P}_\mathrm{OIRT} \nonumber\\
& & \qquad + \eta_0 \mathbf{Q}_0 + \eta_\mathrm{B}
\mathbf{Q}_\mathrm{B} + \eta_\mathrm{BP} \mathbf{Q}_\mathrm{BP} +
\eta_\mathrm{BPS} \mathbf{Q}_\mathrm{BPS}. \nonumber
\end{eqnarray}
Randomized-inclusive randomizations are used in this experiment, as
the outcome of the randomization of lines to samples must 
be known before the samples can be randomized to analyses.
The Plots pseudofactors~$\mathrm{P}_1$ and~$\mathrm{P}_2$ are 
used to ensure appropriate partial confounding of sources from the
lines tier with sources in the analyses tier.
These pseudofactors do not give idempotents in $\mathbf{V}$, 
because they do not contribute to the variance matrix;
they are irrelevant to the randomization of lines to samples,
and are not among the unrandomized factors, that give rise to
covariance, in the randomization of samples to analyses.
However, as in Example~\ref{eg:cotton},
$\mathbf{Q}_\mathrm{BP}$ can be rewritten as the sum of three
$\mathbf{Q}^*$-matrices each with coefficient~$\eta_{\mathrm{BP}}$.
This results in the coefficient $\eta_{\mathrm{BP}}$ occurring with
three different $\xi$-coefficients in the skeleton
anova in Table~\ref{tab:ANOVAWheat2}, which is an extended version of the
decomposition table given for Example~5 in \citep{Brien092}.
To obtain structure balance, $\image(\mathbf{Q}_\mathrm{BPS})$
is decomposed as the sum of
five subspaces involving the pseudofactor $\mathrm{S}_1$ so that
$\mathbf{Q}_\mathrm{BPS}$
can be rewritten as the sum of five
$\mathbf{Q}^*$-matrices each with coefficient $\eta_{\mathrm{BPS}}$. As
a consequence, the coefficient $\eta_{\mathrm{BPS}}$ occurs with five
different $\xi$-coefficients in Table~\ref{tab:ANOVAWheat2}.
\end{egg}

\begin{egg}[Small example]
\label{eg:artif} Figure~\ref{fig:artif} is the randomization diagram
for the first example in \citep{Wood88}.
The levels of pseudofactor $\mathrm{U}_1$ give the treatments
allocated to units in the first phase.  The small open circle
indicates the nonorthogonal block design (group-divisible for
$\mathrm{U}_1$) for the second phase. As the design allocating units
to plots in the second phase depends on the outcome of the
randomization of treatments to units in the first phase, 
the randomizations are randomized-inclusive.

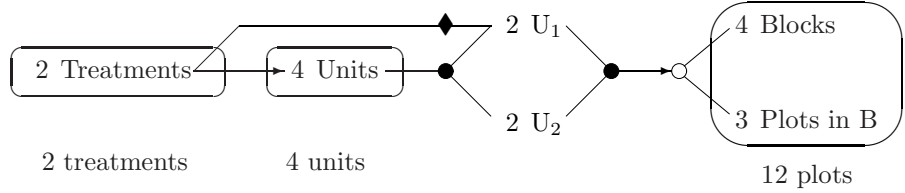
\begin{figure}[htbp]
\setlength{\unitlength}{1cm}
\small
\centering
\begin{picture}(13.8,2.25)(-1,0.2)
\put(0,1.5){\begin{tierbox}2 & Treatments \hspace{250pt}\end{tierbox}}
\put(2.55,1.6){\vector(1,0){1.2}}
\put(2.55,1.6){\line(1,1){0.6}}
\put(3.4,1.5){\begin{tierbox}4 & Units\end{tierbox}}
\put(5.9,1.6){\line(-1,0){0.8}}
\put(5.9,1.6){\blob}
\put(5.9,1.6){\line(1,1){0.6}}
\put(5.9,1.6){\line(1,-1){0.6}}
\put(6.7,2.2){\makebox(1,0)[l]{$2 \hspace{\nlevnamesep} {\rm U}_1$}}
\put(6.7,0.9){\makebox(1,0)[l]{$2 \hspace{\nlevnamesep} {\rm U}_2$}}
\put(3.15,2.2){\line(1,0){3.35}}
\put(5.78,2.12){\makepseudo}
\put(8.1,1.6){\line(-1,1){0.6}}
\put(8.1,1.6){\line(-1,-1){0.6}}
\put(8.1,1.6){\blob}
\put(8.1,1.6){\vector(1,0){0.8}}
\put(9.0,1.6){\nonorthcircle}
\put(9.07,1.67){\line(6,5){0.6}}
\put(9.07,1.53){\line(6,-5){0.6}}
\put(9.3,1.475){\begin{tierbox}4 & Blocks\\ \\ \\ 3 & Plots in B\end{tierbox}}
\put(1.5,0.4){\makebox(0,0){2 treatments}}
\put(4.3,0.4){\makebox(0,0){4 units}}
\put(10.7,0.2){\makebox(0,0){12 plots}}
\end{picture}
\caption{\label{fig:artif}Randomization diagram for  Example~\ref{eg:artif}: treatments are randomized to units, which are in
  turn randomized to plots; $\mathrm{B}$ denotes Blocks;
  $\mathrm{U}_1$ and $\mathrm{U}_2$ are pseudofactors for Units,
  forming groups of units 
determined by the treatments.}
\end{figure}

The skeleton anova in Table~\ref{tab:artif} agrees with
the conclusions reported in \citep{Wood88}.  However, it is more
informative. For example, using Section~\ref{s:single},
the variance of the within-blocks estimator of the treatment
difference is equal to 
$(2/6) \times (\xi_{\mathrm{BP}} +(8/3)\eta_{\mathrm{U}}) /(8/9)$.
Table~\ref{tab:artif} shows that an estimator of $\xi_{\mathrm{BP}}
+(8/3)\eta_{\mathrm{U}}$ is twice the mean square for
${\mathrm{Plots}\nesting{\mathrm{B}}} \combine (\mathrm{U} \resid
\mathrm{U}_1)$ minus the Residual mean square. 

\begin{table}[h]
\caption{\label{tab:artif}Skeleton analysis of variance
for Example~\ref{eg:artif}}
\begin{center}
\begin{tabular}{lrc*{2}{clrc}l}
\toprule
\multicolumn{2}{c}{\textbf{plots tier}} & &
\multicolumn{3}{c}{\textbf{units tier}} & &
\multicolumn{3}{c}{\textbf{treatments tier}} & \\
\cmidrule{1-2} \cmidrule{4-6} \cmidrule{8-10}
\textbf{source} & \textbf{d.f.} & &
\textbf{eff.} & \textbf{source} & \textbf{d.f.} & &
\textbf{eff.} & \textbf{source} & \textbf{d.f.} & &
\multicolumn{1}{c}{\textbf{E.M.S.}} \\
\midrule
Mean & 1 & & & Mean & 1 & & & Mean & 1 & & $\xi_0 + 3\eta_0 + q_0$ \\
\midrule
Blocks & 3 & & \raisebox{8.5pt}[9pt][5pt]{}$\frac{1}{9}$ &
    $\mathrm{U}_1$ & 1 & &
\raisebox{8.5pt}[9pt][5pt]{}$\frac{1}{9}$ & Treatments & 1 & &
               $\xi_\mathrm{B} +
               \raisebox{8.5pt}[9pt][5pt]{}\frac{1}{9}3\eta_\mathrm{U}
               + \raisebox{8.5pt}[9pt][5pt]{}\frac{1}{9}q(\mathrm{T})$ \\
 & & & \raisebox{8.5pt}[9pt][5pt]{}$\frac{5}{9}$ &
    $\mathrm{U} \resid \mathrm{U}_1$ & 2 & & &  & & &
               $\xi_\mathrm{B} +
               \raisebox{8.5pt}[9pt][5pt]{}\frac{5}{9}3\eta_\mathrm{U}$  \\
\midrule
$\mathrm{Plots} \nesting{\mathrm{B}}$ & 8 & &
\raisebox{8.5pt}[9pt][5pt]{}$\frac{8}{9}$ & $\mathrm{U}_1$ & 1 & &
\raisebox{8.5pt}[9pt][5pt]{}$\frac{8}{9}$ & Treatments & 1 & &
               $\xi_\mathrm{BP} +
               \raisebox{8.5pt}[9pt][5pt]{}\frac{8}{9}3\eta_\mathrm{U} +
               \raisebox{8.5pt}[9pt][5pt]{}\frac{8}{9}q(\mathrm{T})$ \\
 & & & \raisebox{8.5pt}[9pt][5pt]{}$\frac{4}{9}$ &
    $\mathrm{U} \resid \mathrm{U}_1$ & 2 & & & & & &
               $\xi_\mathrm{BP} +
               \raisebox{8.5pt}[9pt][5pt]{}\frac{4}{9}3\eta_\mathrm{U}$  \\
& & & & Residual  & 5 & & & & & & $\xi_\mathrm{BP} $  \\
\bottomrule
\end{tabular}
\end{center}
\end{table}
\end{egg}

\section{Estimation in a two-tiered experiment}
\label{sec:esti2}

Estimation of treatment effects and variances is straightforward in a
two-tiered experiment with structure balance.

\subsection{Estimating treatment effects and variances in one stratum}
\label{sec:orthcase}

For data satisfying the conditions in Section~\ref{sec:aov},
the following are also shown in
\citep{Bailey81,Houtman83,Nelder65b}.
\begin{enumerate}
\item[(E.a)]
The best linear unbiased estimator of the 
  treatment effects
  $\mathbf{R}\mathbf{X}_h \boldsymbol{\tau}$, using only the
  projected data $\mathbf{QY}$, is $\mathbf{R}(\mathbf{Q}
  \combine \mathbf{R}) \mathbf{Y} /\lambda_{\mathbf{QR}}$,
  which is equal to $\mathbf{R}\mathbf{Q} \mathbf{Y} /\lambda_{\mathbf{QR}}$.
\item[(E.b)]
The variance matrix of the above estimator is
  $(\eta_\mathbf{Q}/\lambda_{\mathbf{Q}\mathbf{R}})\mathbf{R}$.
\item[(E.c)]
From (A.d) in Section~\ref{sec:aov}, an unbiased
estimate of $\eta_\mathbf{Q}$ is given by the mean square
for $\mathbf{Q} \resid \mathcal{R}$, if $\mathbf{Q} \resid \mathcal{R}$ is 
nonzero.
\end{enumerate}

\subsection{Treatment structure orthogonal to variance structure}
\label{sec:allorth}

If $\mathcal{R}$ is orthogonal in relation to $\mathcal{Q}$ then
each $\mathbf{R}$ in $\mathcal{R}$ has
some $\mathbf{Q}$ in $\mathcal{Q}$ such that $\lambda_{\mathbf{QR}}=1$.
Then all the information on $\mathbf{R} \mathbf{X}_h \boldsymbol{\tau}$
is in stratum $\image (\mathbf{Q})$.
Hence result~(E.a) in Section~\ref{sec:orthcase} gives $\mathbf{RQY}$ as the
overall best linear unbiased estimator of
$\mathbf{R} \mathbf{X}_h \boldsymbol{\tau}$.
Result~(E.b) shows that the variance matrix
of this estimator is $\eta_{\mathbf{Q}}\mathbf{R}$, and
result~(E.c) that the mean square for $\mathbf{Q} \resid \mathcal{R}$
is an unbiased estimator for $\eta_{\mathbf{Q}}$, if $\mathbf{Q}
\resid \mathcal{R}$ is nonzero.

\subsection{Estimating treatment effects from multiple strata when
variances are known}
\label{sec:combine}

Suppose that $\mathcal{R}$ is not orthogonal in relation to $\mathcal{Q}$.
As shown in \citep{Houtman83,Nelder68},
if the coefficients
$\eta_{\mathbf{Q}}$ are known then we can combine information on
$\mathbf{R}\mathbf{X}_h \boldsymbol{\tau}$ from all strata for which
$\lambda_{\mathbf{QR}} \ne 0$ to obtain its generalized least squares
(GLS) estimator, which is the best linear unbiased estimator.
In our notation, it is given by
\begin{equation}\label{eq:GLS}
\mathbf{R}\mathbf{X}_h \widehat{\boldsymbol{\tau}} =
\theta^{-1}_{\mathbf{R}}
\sum_{\mathbf{Q} \in \mathcal{Q}}
\eta^{-1}_{\mathbf{Q}} \mathbf{R}\mathbf{Q} \mathbf{Y},
\end{equation}
where $\theta_{\mathbf{R}} =
\sum_{\mathbf{Q} \in \mathcal{Q}}
\lambda_{\mathbf{QR}} \eta^{-1}_{\mathbf{Q}}$.  The
variance matrix of this estimator is
$\theta^{-1}_{\mathbf{R}}\mathbf{R}$.

\subsection{Estimating treatment effects and variances from multiple strata}
\label{sec:estall}

However, usually the coefficients $\eta_{\mathbf{Q}}$ are unknown and
must be estimated.
One method for this is to use the mean square for $\mathbf{Q} \resid
\mathcal{R}$ to estimate $\eta_{\mathbf{Q}}$.
Nelder \citep{Nelder68} argued that, especially for
designs in which some strata have few Residual degrees of freedom,
estimates should instead be obtained by equating
the expected and observed values of the mean squares
for what Houtman and Speed \citep{Houtman83} called `actual 
residuals'.
Even though normality is not assumed, these estimates are the same as those
obtained by REML \cite{Patterson71}
because, as
is shown in \citep{Houtman83}, Section~4.5, 
and \citep{Patterson75}, the same
set of equations has to be solved for both.
As will be discussed in Section~\ref{sec:mixedrando}, the constraints
on the variance parameters being estimated here are 
different from those for a variance-components model.

As noted in \citep{Houtman83,Nelder68},
the estimation of the coefficients $\eta_{\mathbf{Q}}$ requires
an iter\-ative procedure, because their
estimation needs the estimated value of $\boldsymbol{\tau}$
and vice versa.
Given working estimates $\widehat{\eta}^*_{\mathbf{Q}}$ 
of $\eta_{\mathbf{Q}}$,
a working estimate $\widehat{\boldsymbol{\tau}}^*$ 
of $\boldsymbol{\tau}$ can be obtained
from equation~\eqref{eq:GLS}: thus a revised estimate
of each $\eta_{\mathbf{Q}}$ can be computed as
\begin{equation}
\label{eqn:xiest}
   \frac{\mathbf{y}' \left( \mathbf{Q} \resid \mathcal{R} \right ) \mathbf{y} +
   \left ( \summ_{\mathbf{R} \in \mathcal{R}}
         \left \{ \mathbf{y}' (\mathbf{Q} \combine \mathbf{R}) \mathbf{y}
             - \lambda_{\mathbf{QR}}
               \widehat{\boldsymbol{\tau}}^{*_{'}} \mathbf{X}_h' \mathbf{R}
               \mathbf{X}_h \widehat{\boldsymbol{\tau}}^*
         \right \}
   \right )}
   {d'_{\mathbf{Q}}},
\end{equation}
where $\summ_{\mathbf{R} \in \mathcal{R}}$ means summation over
$\mathbf{R} \in \mathcal{R}$ for which $\lambda_{\mathbf{QR}} \neq 0$, and
$d'_{\mathbf{Q}}$ are 
the \emph{effective} degrees of freedom for this estimator,
which are given by
\[
d'_{\mathbf{Q}} =
  \trace{(\mathbf{Q} \resid \mathcal{R})} + \summ_{\mathbf{R} \in \mathcal{R}}
  \left[ 1 - \theta^{-1}_{\mathbf{R}} (\widehat{\eta}^*_{\mathbf{Q}})^{-1}
  \lambda_{\mathbf{QR}} \right ]
  \trace{\mathbf{R}} .
\]
Since $\mathbf{Q} \resid \mathcal{R}$ and $\mathbf{R}$
are both idempotent, their traces are equal to their ranks.

The numerator of expression~\eqref{eqn:xiest} is the sum of two parts.
The first is the Residual sum of squares in this stratum from the
anova; the second is the difference between
the sum of squares of the treatment estimates from just the data
projected onto $\image(\mathbf{Q})$ and the sum of squares of
the combined estimates, summed over all $\mathbf{R}$ for which
$\lambda_{\mathbf{QR}} \neq 0$.
The former does not depend on $\eta_{\mathbf{Q}}$, but the latter does.
The effective degrees of freedom make it clear that, even when
$\mathbf{Q} \resid \mathcal{R} = \mathbf{0}$, there can be
information to estimate $\eta_{\mathbf{Q}}$.

If estimates of the canonical components are required these can be
obtained from the estimates of the spectral components.


\section{Estimating treatment effects and variances in a single
part of ${\mathcal{P} \combine \mathcal{Q}}$}
\label{s:single}

Suppose that $\lambda_{\mathbf{PQ}} \ne 0$, so that there is an
idempotent $\mathbf{P} \combine \mathbf{Q}$.  Consider an
idempotent~$\mathbf{R}$ in $\mathcal{R}$ for which
$\lambda_{\mathbf{P}\combine \mathbf{Q}, \mathbf{R}} \ne 0$.
Theorem~5.1 of \citep{Brien092} shows that
$\lambda_{\mathbf{P}\combine \mathbf{Q}, \mathbf{R}} =
\lambda_{\mathbf{PQ}} \lambda_{\mathbf{QR}}$.
Applying the results of Section~\ref{sec:orthcase} with $\mathbf{Y}$
 and $\mathbf{Q}$ replaced by $\mathbf{PY}$ and $\mathbf{P} \combine
 \mathbf{Q}$ respectively, and using equation \eqref{eq:pdata} for
 $\Cov{\mathbf{PY}}$, we find that the best linear
unbiased estimator of the treatment effect
$\mathbf{R}\mathbf{X}_f \mathbf{X}_h \boldsymbol{\tau}$, using only
the projected data $(\mathbf{P} \combine \mathbf{Q})\mathbf{Y}$, is
$\mathbf{R}((\mathbf{P} \combine \mathbf{Q})\combine \mathbf{R})
\mathbf{Y} /\lambda_{\mathbf{PQ}}\lambda_{\mathbf{QR}}$, which is
equal to $\mathbf{R} \mathbf{Q} \mathbf{P} \mathbf{Y}/
\lambda_{\mathbf{PQ}} \lambda_{\mathbf{QR}}$.
Moreover, the variance matrix of this estimator is equal to
$\mathbf{R}
(\xi_\mathbf{P} + r\lambda_{\mathbf{PQ}}\eta_\mathbf{Q})/
\lambda_{\mathbf{PQ}}\lambda_{\mathbf{QR}}$.
Result (A.i) in Section~\ref{sec:aov3} shows that the mean square for
$(\mathbf{P} \combine \mathbf{Q}) \resid \mathcal{R}$ is an unbiased
estimator for $\xi_\mathbf{P} + r\lambda_{\mathbf{PQ}}\eta_\mathbf{Q}$,
if $(\mathbf{P} \combine \mathbf{Q}) \resid \mathcal{R}$ is nonzero.

In Example~\ref{eg:2PhaseSensory}, the effects for M and T$\inter$M
are estimated in just the source $\mathrm{P} \nesting{\mathrm{O}
  \wedge \mathrm{I} \wedge \mathrm{S} 
\wedge \mathrm{J}}\combine \mathrm{H} \nesting{\mathrm{Q} \wedge
\mathrm{R} \wedge \mathrm{C}}$ and the Residual mean square for
$\mathrm{H} \nesting{\mathrm{Q} \wedge \mathrm{R} \wedge \mathrm{C}}$
is an unbiased estimator of $\xi_\mathrm{OISJP} + 12\eta_\mathrm{QRCH}$.

\section{Full estimation in a three-tiered experiment which is
anova-applicable}
\label{sec:good}

\subsection{Full or partial anova}
\label{sec:fullpart}
Call the triple $(\mathcal{P}, \mathcal{Q}, \mathcal{R})$
\emph{anova-applicable}
if it satisfies the following condition:
\begin{equation}
\mbox{for every $\mathbf{Q}$ in $\mathcal{Q}$,
if $\mathbf{Q} \mathbf{I}_{\mathcal{R}} \ne \mathbf{0}$
then $\mathbf{Q}\in\mathcal{Q}_1$.}
\label{anovaapp}
\end{equation}
That is, if the source for an $\mathbf{R}$ in $\mathcal{R}$ is
(partially) confounded with a source for some $\mathbf{Q}$, then 
the latter source must be confounded with the source corresponding to a
single $\mathbf{P}$.
Section~\ref{sec:aov3} shows that when this condition is satisfied
then the idempotents in
$\mathcal{P} \combine \mathcal{Q}$ whose subspaces are 
contained in eigenspaces of $\mathbf{V}$  include all those which
have any part of $\mathcal{R}$ partially or totally confounded with them.

Condition~\eqref{anovaapp} is satisfied when $\mathcal{Q}$ is orthogonal in
relation to $\mathcal{P}$, so that $\mathcal{Q}_1 = \mathcal{Q}$. Then
$\mathbf{V}$ is  given by equation~\eqref{eq:OVS}, possibly with OVS.
Estimation of treatment effects and their variances proceeds as in
Section~\ref{sec:esti2}.
Examples~\ref{eg:tbb}, \ref{eg:cotton} and~\ref{eg:Wheat2} are like
this. We call this \emph{full anova}.

Under full anova, if no
$(\mathbf{P} \combine \mathbf{Q}) \resid \mathcal{R} = \mathbf{0}$, we
estimate linear combinations of spectral components from the anova,
even if they are not needed for standard errors of treatment effects.
Otherwise, formula~\eqref{eqn:xiest} 
can be used, but with $\mathbf{P} \combine \mathbf{Q}$
replacing $\mathbf{Q}$.
In Example~\ref{eg:tbb}, all  the information about each treatment
source in $\mathcal{R} \setminus \mathbf{R}_0$ is contained in
$(\mathrm{P}\inter\mathrm{T}\nesting{\mathrm{S}})
  \combine \mathrm{M}\nesting{\mathrm{B}}$. Also,
the difference between the  mean squares for
$((\mathrm{P}\inter\mathrm{T}\nesting{\mathrm{S}})
    \combine (\mathrm{M}\nesting{\mathrm{B}}))\resid \mathcal{R}$ and
$(\mathrm{P}\inter\mathrm{T}\nesting{\mathrm{S}}) \resid \mathcal{Q}$
estimates $12\eta_\mathrm{BM}$.

In general, put $\mathcal{P} * \mathcal{Q} =
\mathcal{Q}_1 \cup \{\mathbf{P}
\resid \mathcal{Q}: \mathbf{P}\in \mathcal{P} \}$.
Then the images of all the idempotents in $\mathcal{P} * \mathcal{Q}$
are contained in eigenspaces of $\mathbf{V}$.
If $(\mathcal{P}, \mathcal{Q}, \mathcal{R})$ is
anova-applicable but $\mathcal{Q}$ is not orthogonal in relation to
$\mathcal{P}$ then we have \emph{partial anova}, using only the
information in $\mathcal{P} * \mathcal{Q}$.
A treatment idempotent $\mathbf{R}$ in~$\mathcal{R}$ may be 
nonorthogonal to more than one part of $\mathcal{P} \combine
\mathcal{Q}$, but these are all in $\mathcal{P} * \mathcal{Q}$.
Section~\ref{sec:aov3} shows that estimators of variances of treatment
effects which are in different parts of $\mathcal{P} * \mathcal{Q}$
are uncorrelated, and so information can be combined as in
Section~\ref{sec:estall}.

However, the linear combinations of
spectral components in the expected mean square for parts of
$\mathcal{P} \combine \mathcal{Q}$ outside $\mathcal{P} * \mathcal{Q}$
are not involved in this process, and their anova-estimators may not
have good properties.
A similar situation arises in two-tiered experiments if the group is
not stratifiable but 
all treatment subspaces are contained in known eigenspaces of the
variance matrix (see Example~16 in \citep{Bailey91}).

In the special case that $\mathcal{R}$ is orthogonal in
relation to $\mathcal{P}\combine\mathcal{Q}$,
each $\mathbf{R}$ in~$\mathcal{R}$ has unique
idempotents $\mathbf{Q}$ in $\mathcal{Q}$ and $\mathbf{P}$ in
$\mathcal{P}$ such that $\lambda_{\mathbf{PQ}}\lambda_{\mathbf{QR}} =1$,
so that $\lambda_{\mathbf{PQ}} = \lambda_{\mathbf{QR}} = 1$.
Hence, Condition~\eqref{anovaapp} is satisfied.
Then  $\mathbf{P} \combine \mathbf{Q} = \mathbf{Q}$, $\mathbf{RQP} =
\mathbf{R}$ and the effect $\mathbf{R}\mathbf{X}_f
\mathbf{X}_h \boldsymbol{\tau}$ is estimated in just
the one part $\mathbf{P} \combine \mathbf{Q}$ of
$\mathcal{P} * \mathcal{Q}$.  The anova may be either full or partial.
The full anova tables obtained by putting data into each of
Tables~\ref{tab:tbb} and~\ref{tab:cox4} are like this.
The estimator for $\mathbf{R}\mathbf{X}_f
\mathbf{X}_h \boldsymbol{\tau}$,  and its variance matrix,
obtained by simplifying the expressions given in Section~\ref{s:single}, are
$\mathbf{RY}$ and $(\xi_{\mathbf{P}} + r\eta_{\mathbf{Q}} )\mathbf{R}$,
respectively. If $\xi_{\mathbf{P}}$ and $\eta_{\mathbf{Q}}$ are not known,
the mean square for $(\mathbf{P} \combine \mathbf{Q}) \resid \mathcal{R}$
provides an unbiased estimate of $\xi_{\mathbf{P}} +
r\eta_{\mathbf{Q}}$,
unless $(\mathbf{P} \combine \mathbf{Q}) \resid \mathcal{R} = \mathbf{0}$.
Estimation of treatment effects and their standard errors proceeds
exactly as in Section~\ref{sec:allorth}.

\subsection{Difficulties that do not arise in two-tiered experiments}

Even when $(\mathcal{P}, \mathcal{Q}, \mathcal{R})$ is
anova-applicable, some phenomena can
occur that are not possible in two-tiered experiments, even for the
straightforward special case where $\mathcal{P}$ and $\mathcal{Q}$ are
both  poset block structures, $\mathcal{Q}$ is orthogonal in relation
to $\mathcal{P}$, and $\mathcal{R}$ is orthogonal in relation to
$\mathcal{Q}$. 

\subsubsection{Inestimability of some spectral and canonical components}
\label{sec:inesti}

For a two-tiered experiment with OVS and $\mathcal{R}$ orthogonal in
relation to $\mathcal{Q}$, 
the estimability of spectral components of variance is easily
determined.  If $\mathbf{Q} \resid \mathcal{R}$ is nonzero then its
mean square provides the best unbiased quadratic estimator of
$\eta_{\mathbf{Q}}$; otherwise, there is no estimator for
$\eta_{\mathbf{Q}}$.  In particular, $\eta_0$ is never estimable.

In a three-tiered experiment, the expected mean square for each
Residual source is a linear combination of a $\xi$-parameter and an
$\eta$ parameter. It may not be possible to estimate $\xi$ and $\eta$
separately.  This affects 
the estimability of canonical components, although it appears that
often more individual canonical, than spectral, components are estimable.
The parameters $\xi_0$, $\eta_0$, $\phi_0$ and $\psi_0$ are never estimable.

Otherwise, the simplest way in which two spectral components cannot be
estimated separately occurs when a generalized 
factor $F$ on $\Upsilon$ is randomly assigned to a generalized
factor $H$ on $\Omega$ with the same number of levels.
Then only a linear combination of $\xi_\mathrm{H}$
and $\eta_{\mathrm{F}}$ can be estimated,
and hence only a linear combination of $\phi_H$ and $\psi_F$.
In Example~\ref{eg:tbb},
Blocks are assigned to Sessions, both of which have three levels.
As a result, only $\xi_\mathrm{S} + 12\eta_\mathrm{B}$ is estimable,
as is shown in Table~\ref{tab:tbb}, where 
these two components only occur together. Correspondingly, 
only $\phi_\mathrm{S} + \psi_\mathrm{B}$ is estimable.

In the special case that $\card{\Upsilon} = \card{\Omega}$,
Lemma~4.2 of \citep{Brien092} shows that
$\mathcal{P}\combine\mathcal{Q} = \mathcal{Q}$ and so
there are no idempotents of the form $\mathbf{P}\resid\mathcal{Q}$. 
Thus every expected mean square contains
one $\eta_{\mathbf{Q}}$ and one $\xi_{\mathbf{P}}$.
If $\kappa$ is any constant smaller in modulus than all the
$\eta_{\mathbf{Q}}$ and all the $\xi_{\mathbf{P}}$, then 
$\kappa\mathbf{I}_{\Omega}$ can be added to
$\sum\xi_{\mathbf{P}} \mathbf{P}$ and subtracted from
$\sum\eta_{\mathbf{Q}} \mathbf{Q}$ without changing the
variance matrix $\mathbf{V}$ in equation~\eqref{eq:twocov}.
Thus none of the spectral components of variance can be estimated,
although sums of the form $\xi_{c(\mathbf{Q})} + \eta_{\mathbf{Q}}$ can be.
For estimates of standard errors,
these sums are all that is needed, and so there is no problem. 
However, for comparing sources of variation,
estimates of canonical components are required. Except for $\phi_{\Omega}$ and 
$\psi_{\Upsilon}$, each canonical component is a 
multiple of a difference between spectral components.
This may well be estimable, even though the corresponding spectral
component is not. 

In Example~\ref{eg:cotton}, none of the spectral components is estimable.
If relative magnitudes of sources of variation are to be investigated, 
then the canonical components are needed:
all of these are estimable except for $\phi_\mathrm{OT}$ and 
$\psi_\mathrm{BPF}$, whose sum is estimable, and $\phi_0$ and $\psi_0$.
Thus the only real restriction in estimating 
canonical components in this experiment is that it is not possible to
separate variability arising from Tests and Fibres. 
Plan~A for Example~4 of \citep{BrBa:mult} is an
alternative design for this experiment.
It has no pseudosources, but exhibits the same inestimability.

\subsubsection{Negative estimates of spectral components}
\label{sec:nonneg}

As noted in Section~\ref{sec:derive},
all spectral components of variance must be nonnegative 
and hence so must any linear combination with positive
coefficients. However, if one mean square
involving both a $\xi$ and an $\eta$ is less than that involving just
the same $\xi$, then the anova-estimate 
of $\eta$ is negative. This is analogous to a negative estimate of a
nonnegative variance component; the component should be set to zero. 

In Example~\ref{eg:tbb}, Table~\ref{tab:tbb} shows that the
appropriate Residual source for all three treatment sources is the one
with expected mean square equal to
$\xi_{\mathrm{STP}} + 12\eta_{\mathrm{BM}}$.
Suppose that this mean square turns out to be smaller
than the one whose expectation is $\xi_{\mathrm{STP}}$. 
Then we set $\eta_{\mathrm{BM}}$ to zero and
combine the two Residual mean squares to obtain a better estimate of
$\xi_{\mathrm{STP}}$. 
See Section~\ref{sec:mixmod} for further discussion. 

\subsubsection{The effect of pseudosources}
\label{sec:pseudo}

If there are pseudosources for $\mathcal{Q}$ then some
$\eta$-coefficients occur with more than one $\xi$-coefficient, even
if $\mathcal{Q}$ is orthogonal in relation to $\mathcal{P}$.
This can lead to what we call \emph{linearly dependent commutative
  variance structure} (LDCVS), in which the
eigenspaces of $\mathbf{V}$ are known but the eigen\-values satisfy some
linear equations. This gives a set of sources 
whose expected mean squares  are linearly dependent:  simply equating
them all to their data mean squares may give inconsistent results.
Suppose that for $i$, $j$ in $\{1,2\}$ the idempotent
  $\mathbf{Q}^*_{ij}$ corresponds to a pseudosource for $\mathbf{Q}_j$
  and is totally confounded with $\mathbf{P}_i$.  Then
the expected mean squares for the four idempotents
$(\mathbf{P}_i \combine \mathbf{Q}^*_{ij}) \resid \mathcal{R}$
are $\xi_1 + r\eta_1$, $\xi_1 + r\eta_2$, $\xi_2 +
r\eta_1$ and $\xi_2 + r\eta_2$.
To estimate either the spectral or canonical components by equating
expected and observed mean squares requires that the  the sum of the
middle two observed mean squares is equal to the sum of the outer two.
This situation is similar to the case, discussed by \citet{Bailey12}, of a
mixed model with LDCVS that is orthogonal to the treatment
structure. They show that all four combinations of the following can
occur: all variance components can be estimable, or not; and each
estimable component may, or may not, have a unique estimator
using the mean squares in the anova.

\section{Estimation in a three-tiered experiment which is not
anova-applicable}
\label{sec:comk}

If $\mathbf{V}$ is known then we can estimate treatment effects by
GLS, which gives different results from ordinary least squares if
$\mathcal{R}$ is not orthogonal in relation to $\mathcal{P} \combine
\mathcal{Q}$. 

Equation~\eqref{eq:twocov} gives
\begin{equation*}
\mathbf{V} = \mathbf{I} \mathbf{V}\mathbf{I} =
\sum_{\mathbf{P}\in\mathcal{P}} \xi_{\mathbf{P}} \mathbf{P} +
r \sum_{\mathbf{P}\in\mathcal{P}} \sum_{\mathbf{Q}\in\mathcal{Q}}
\sum_{\mathbf{P^*}\in\mathcal{P}} \eta_{\mathbf{Q}}
\mathbf{P}\mathbf{Q}\mathbf{P}^*.
\end{equation*}
Put $\alpha_{\mathbf{Q}} = \sum_{\mathbf{P}} (\lambda_{\mathbf{PQ}}/\xi_{\mathbf{P}})$. Then
direct calculation shows that
\[
\mathbf{V}^{-1} = \sum_{\mathbf{P}\in\mathcal{P}}
\frac{1}{\xi_{\mathbf{P}}} \mathbf{P}
 - \sum_{\mathbf{P}\in\mathcal{P}} \sum_{\mathbf{Q}\in\mathcal{Q}}
\sum_{\mathbf{P^*}\in\mathcal{P}}
\left( \frac{r\eta_{\mathbf{Q}}}{1+
  r\eta_{\mathbf{Q}}\alpha_{\mathbf{Q}}} \right)
 \frac{1}{\xi_{\mathbf{P}}\xi_{\mathbf{P}^*}} \mathbf{P}\mathbf{Q}\mathbf{P}^*.
\]

Consider $\mathbf{R}$ in $\mathcal{R}$.
When the $\xi_{\mathbf{P}}$ and $\eta_{\mathbf{Q}}$ are known, the GLS 
estimator of the treatment effect
$\mathbf{R}\mathbf{X}_f \mathbf{X}_h\boldsymbol{\tau}$
is $(\mathbf{R}\mathbf{V}^{-1}\mathbf{R})^{-}
\mathbf{R} \mathbf{V}^{-1} \mathbf{Y}$,  with
variance matrix $(\mathbf{R} \mathbf{V}^{-1}
\mathbf{R})^{-}$.  For a chain of randomizations,
$\mathbf{R} =\mathbf{R} \mathbf{I}_{\cal Q} =
\sum_{\mathbf{Q}} \mathbf{R} \mathbf{Q}$, so
\begin{eqnarray*}
\mathbf{R} \mathbf{V}^{-1}  &=&
\sum_{\mathbf{P}} \sum_{\mathbf{Q}} \frac{1}{\xi_{\mathbf{P}}}
\mathbf{R} \mathbf{Q}\mathbf{P}  -
\sum_{\mathbf{P}} \sum_{\mathbf{Q}} \sum_{\mathbf{P}^*}
\left( \frac{r\eta_{\mathbf{Q}}}{1+
  r\eta_{\mathbf{Q}}\alpha_{\mathbf{Q}}} \right)
 \frac{1}{\xi_{\mathbf{P}}\xi_{\mathbf{P}^*}}
 \mathbf{R}\mathbf{Q}\mathbf{P}\mathbf{Q}\mathbf{P}^*\\
  &=&
\sum_{\mathbf{P}} \sum_{\mathbf{Q}} \frac{1}{\xi_{\mathbf{P}}}
\mathbf{R} \mathbf{Q}\mathbf{P}  -
 \sum_{\mathbf{Q}} \sum_{\mathbf{P}^*}
\left( \frac{r\eta_{\mathbf{Q}}\alpha_{\mathbf{Q}}}{1+
  r\eta_{\mathbf{Q}}\alpha_{\mathbf{Q}}} \right)
 \frac{1}{\xi_{\mathbf{P}^*}}
 \mathbf{R}\mathbf{Q}\mathbf{P}^*\\
&=& \sum_{\mathbf{P}} \sum_{\mathbf{Q}}
 \frac{1}{\xi_{\mathbf{P}}}
\left(\frac{1}{1+ r\eta_{\mathbf{Q}}\alpha_{\mathbf{Q}}} \right)
\mathbf{R} \mathbf{Q}\mathbf{P}.
\end{eqnarray*}
Hence
\[
\mathbf{R} \mathbf{V}^{-1} \mathbf{R} = 
\sum_{\mathbf{P}} \sum_{\mathbf{Q}}
 \frac{1}{\xi_{\mathbf{P}}}
\left(\frac{1}{1+ r\eta_{\mathbf{Q}}\alpha_{\mathbf{Q}}} \right)
\mathbf{R} \mathbf{Q}\mathbf{P} \mathbf{Q} \mathbf{R} =
\sum_{\mathbf{Q}} \left(\frac{\alpha_{\mathbf{Q}} \lambda_{\mathbf{QR}}}
{1+ r\eta_{\mathbf{Q}}\alpha_{\mathbf{Q}}} \right)  \mathbf{R} 
 = \theta_{\mathbf{R}} \mathbf{R}, 
\]
with $\theta_{\mathbf{R}} = \sum_{\mathbf{Q}\in\mathcal{Q}}
\alpha_{\mathbf{Q}}\lambda_{\mathbf{QR}}(1 + r\eta_{\mathbf{Q}}\alpha_{\mathbf{Q}})^{-1}$.
Thus the GLS  estimator of
$\mathbf{R}\mathbf{X}_f \mathbf{X}_h \boldsymbol{\tau}$ is
\[
\frac{1}{\theta_{\mathbf{R}}}
\sum_{\mathbf{P}} \sum_{\mathbf{Q}}
 \frac{1}{\xi_{\mathbf{P}}}
\left(\frac{1}{1+ r\eta_{\mathbf{Q}}\alpha_{\mathbf{Q}}} \right)
\mathbf{R} \mathbf{Q}\mathbf{P} \mathbf{Y},
\]
with variance matrix $\theta_{\mathbf{R}}^{-1} \mathbf{R}$.

This estimator is a linear combination of the $\mathbf{R}\mathbf{Q}
\mathbf{P} \mathbf{Y}$.
There are no terms in $\mathbf{R}(\mathbf{P} \resid \mathcal{Q})$,
because all projectors of this form are zero for a chain of randomizations.

In the special case that $\cal R$ is orthogonal in
relation to $\cal Q$ 
there is a unique~$\mathbf{Q}$ such that
$\mathbf{R}\mathbf{Q} = \mathbf{R}$ while $\mathbf{R} \mathbf{Q}^* =
\mathbf{0}$ if $\mathbf{Q}^*\ne \mathbf{Q}$, so the estimator is a
linear combination of the
$\mathbf{R}\mathbf{P}\mathbf{Y}$, as shown in \citep{Wood88}.
The scalar $\theta_{\mathbf{R}}$ 
specializes to that given  in \citep{Wood88}.

In the anova-applicable case, we have
$\theta_{\mathbf{R}} = \sum_{\mathbf{Q} \in \mathcal{Q}}
\lambda_{\mathbf{QR}}(\xi_{c(\mathbf{Q})} + r\eta_{\mathbf{Q}})^{-1}$ and
\begin{equation*}
\mathbf{R} \mathbf{X}_f \mathbf{X}_h \widehat{\boldsymbol{\tau}} =
\theta^{-1}_{\mathbf{R}} \sum_{\mathbf{Q}\in\mathcal{Q}}
 (\xi_{c(\mathbf{Q})} + r\eta_{\mathbf{Q}})^{-1}
\mathbf{R} \mathbf{Q} \mathbf{Y}.
\end{equation*}

If $\mathbf{V}$ is unknown, canonical components need to be estimated
using REML followed by estimation of the 
treatment effects using estimated GLS (EGLS).

\section{Extension to more than two randomizations in a chain}
\label{sec:multi}

We have seen that, even with structure balance, there can be
difficulties with anova estimation for a three-tiered experiment.  One
solution can be to use software designed for fitting mixed models.
This does not need to be restricted to designs with structure balance,
or to experiments with three tiers, so we begin by generalizing
Section~\ref{sec:ranmod}. 

As discussed in \citep{BrBa:mult}, Section 6, and \cite{Brien092},
Section~7, more than two randomizations are possible. For example, a
multiphase experiment can consist of $p$ phases and involve $p$
randomizations. Then there are $p$ sets, $\Omega_i$ for $i = 1$,
\ldots, $p$, and another set $\Gamma$ for the first phase. There is a
design function $h\colon\Omega_p\to\Gamma$;  if the objects in 
$\Gamma$ are treatments then $h(\omega)$ is the treatment
assigned to unit $\omega$ in $\Omega_p$.  A stratifiable group $G_p$
of permutations of $\Omega_p$ is used to randomize~$h$. 
For $i=1$, \ldots, $p-1$,
there is a design function $f_i\colon\Omega_{i}\to\Omega_{i+1}$,
so that $f_i(\omega)$ is the unit in $\Omega_{i+1}$ assigned to unit $\omega$
in $\Omega_{i}$; there is also a stratifiable group $G_i$ of permutations
of $\Omega_i$ which is used to randomize $f_i$.
It is assumed that $f_i$ is equireplicate, with replication
$r_{i+1}/r_{i}$, where $r_1=1$, so that each element in $\Omega_i$ is
assigned to $r_i$ elements in $\Omega_1$. 

Let $Y_{\omega}$ be the response on unit~$\omega$ in~$\Omega_1$.
For $\omega$ in $\Omega_1$, put $s_1(\omega) = \omega$, 
$s_{i+1}(\omega) = f_i(s_i(\omega))$ for $i=1$, \ldots, $p-1$
and $t(\omega) = h(s_p(\omega))$.
The randomization-based model in equation~(\ref{eq:2rand}) can be
generalized to 
\[
Y_{\omega} = \sum_{i = 1}^p Z_{i,{s_i(\omega)}} + \tau_{t(\omega)},
\]
where $Z_{i,{s_i(\omega)}}$ is the random effect, under randomization
by $G_i$, for unit~$s_i(\omega)$ in~$\Omega_i$.

For this model,
$\E{\mathbf{Y}} = \mathbf{X}_s\mathbf{X}_h\boldsymbol{\tau}$, where
$\mathbf{X}_s$ is the $\Omega_1 \times \Omega_p$ design matrix for~$s_p$
and $\mathbf{X}_h$ is the $\Omega_p \times \Gamma$ design matrix for~$h$.
Generalize $\mathbf{I}_{\mathcal{R}}$ to be the $\Omega_1 \times
\Omega_1$ matrix of orthogonal projection onto
$\image(\mathbf{X}_s\mathbf{X}_h)$.  Also,
$\mathbf{V} = \sum_{i=1}^p \mathbf{V}_i$, where
$\mathbf{V}_i = r_i \sum_{\mathbf{P}_{ij}\in\mathcal{P}_i}
\xi_{ij} \mathbf{P}_{ij}$
and each $\mathbf{P}_{ij}$ is an idempotent of $\mathbf{V}_i$
with spectral component $\xi_{ij}$.

To this point, there is no need for structure balance, nor do any of
the structures need to be defined by factors.  However, 
if the randomization of $f_i$ is based on a tier of factors
$\mathcal{H}_i$ defining a poset block structure on $\Omega_i$ 
then
$\mathbf{V}_i = \sum_{H \in \mathcal{H}_i} \phi_{H} \mathbf{S}_{H}$, where
$\mathbf{S}_{H}$ is the $\Omega_1 \times \Omega_1$ relationship matrix for
$H$ considered as a factor on $\Omega_1$.

If all of $f_1$, \ldots, $f_p$ and $h$ are structure-balanced then
the results of Sections~\ref{sec:treatsb}, \ref{sec:sameaov},
\ref{s:single}  and \ref{sec:good} can be extended to more than two
randomizations.  In particular, generalize $\mathcal{Q}_1$ be the set of
idempotents $\mathbf{Q}$ in $\mathcal{P}_p$ for which there is an idempotent
$\mathbf{P}_{c_i(\mathbf{Q})}$ in $\mathcal{P}_i$ for $i=1$, \ldots,
$p$ such that $\mathbf{P}_{c_p(\mathbf{Q})}\mathbf{Q} = \mathbf{Q}$ and 
$\mathbf{P}_{c_i(\mathbf{Q})}\mathbf{P}_{c_{i+1}(\mathbf{Q})} =
\mathbf{P}_{c_{i+1}(\mathbf{Q})}$ for $i=1$, \ldots, $p-1$.
The condition for anova-applicability becomes
\[
\mbox{for every $\mathbf{Q}$ in $\mathcal{P}_p$,
if $\mathbf{Q} \mathbf{I}_{\mathcal{R}} \ne \mathbf{0}$
then
$\mathbf{Q}\in\mathcal{Q}_1$.}
\]

\section{Obtaining estimates from data for experiments with a
chain of randomizations}
\label{sec:dataest}

How can standard software be used to obtain, from data, 
estimates of treatment effects and their standard errors and/or 
estimates of canonical components, under randomization-based models? 
Assume that, for $i=1$, \ldots, $p$,
$\mathcal{P}_i$ is given by a poset block structure defined by a
set $\mathcal{H}_i$ of generalized 
factors on $\Omega_i$, which are then expressed as factors on $\Omega_1$.
Two possible procedures,
based on mixed models, are anova and mixed-model fitting.

\subsection{Analysis of variance}
\label{sec:aovanal}

This is the method of choice for  anova-applicable cases in which
the structure $\mathcal{R}$ on $\Gamma$ is also 
orthogonal in relation to $\mathcal{Q}_1$, 
and other cases in which it has been decided that each treatment
effect is to be estimated from a single source, as might be done in
Example~\ref{eg:2PhaseSensory}.
Other anova-applicable cases can be dealt with by anova followed by
combination of information, as in Section~\ref{sec:estall}.

Anova can also be used to estimate canonical components when there 
is CVS and $\mathcal{R}$ is orthogonal in relation to
$\mathcal{Q}$. If there is LDCVS then a 
generalized linear model (GLM) estimates the components.
One fits a GLM to the observed mean squares involved in the estimation.
The GLM has a gamma distribution, identity link, dispersion
parameter equal to $2$, weights equal to the degrees of freedom and
an $\mathbf{X}$ matrix that contains, in each row, the coefficients of
the canonical components for the expected mean square corresponding to
the observed mean square.

The advantage of anova is that it is a noniterative procedure in which
all the quantities are well-defined. Further, nonnegativity
constraints are easily implemented as a manual procedure applied after the
anova has been obtained,
and the inestimability of variance parameters is often inconsequential.

However, most anova software does not produce combined estimates of fixed
effects,  so that it is simpler to use mixed-model fitting.
A further difficulty with anova for multitiered experiments
is that software for it is not generally available, GenStat being the
only package that has specific facilities
\citep{Brien12AMT}. However, it may be possible to specify a set of
terms that will produce the correct decomposition by omitting
some sources. For example, the correct decomposition
is obtained for Example~\ref{eg:cotton} from an anova
or a regression model with the sources Operatives, Blocks, Treatments,
$\mathrm{Plots}\nesting{\mathrm{Blocks}}$ and Residual. This is akin to
fitting a mixed model of convenience, as described in \citep{Brien091},
because it does not contain terms for all the potential sources of
variation that have been identified for the experiment.

\subsection{Classes of mixed models for structures defined by factors}
\label{sec:mixedrando}

Most mixed-model software uses a conditional model
(\citep{Littel06}, Appendix A1, and \citep{Searle92}, Section 4.6):
\[
\mathbf{Y} = \mathbf{X}\boldsymbol{\tau} + \mathbf{ZU} + \mathbf{E},
\]
with $\E{\mathbf{Y} \mid \mathbf{U}} =
\mathbf{X}\boldsymbol{\tau} + \mathbf{ZU}$,
$\Cov{\mathbf{U}} = \mathbf{ZGZ'}$  and
$\Cov{\mathbf{Y} \mid \mathbf{U}} = \Cov{\mathbf{E}}
= \mathbf{R}$,
where $\boldsymbol{\tau}$ is the vector of fixed-effects parameters,
$\mathbf{X}$ is an indicator-variable matrix for fixed effects with
one row for each observation and a column for each fixed effect,
$\mathbf{Z}$ is an indicator-variable matrix with a row for each observation
and a column for each random effect, $\mathbf{U}$ is the vector of
random effects, $\mathbf{E}$ is the vector of random unit effects, and
$\mathbf{G}$ and $\mathbf{R}$ are symmetric matrix functions of the
variance parameters. 
This usage of $\mathbf{R}$ is unrelated to its usage elsewhere in the paper.

This conditional model can be re-expressed in the following marginal form:
\begin{equation}
\label{eq:margmixmodel}
\E{\mathbf{Y}} = \mathbf{X} \boldsymbol{\tau} \mbox{ and }
\Cov{\mathbf{Y}} = \mathbf{V} = \mathbf{ZGZ}' + \mathbf{R}.
\end{equation}
The model for the variance matrix in equation~\eqref{eq:margmixmodel} is referred to as the unstructured variance model; the only condition
imposed is that $\mathbf{V}$ is p.s.d. 

We are concerned with models for the $\Omega_1 \times \Omega_1$
variance matrix $\mathbf{V}$ that are based on sets $\mathcal{H}_i$ of
generalized factors on $\Omega_1$.
Put $\mathcal{H} = \bigcup_{i=1}^p \mathcal{H}_i$.
Let $\mathcal{E}$ be the set of those generalized factors in $\mathcal{H}$
that uniquely index the units in $\Omega_1$,
and let $\mathcal{U} = \mathcal{H} \setminus \mathcal{E}$.
This allows us to write
$\Cov{\mathbf{U}}  = \mathbf{ZGZ'} =
\sum_{U \in \mathcal{U}} \mathbf{Z}_U\mathbf{G}_U\mathbf{Z'}_U$ and
$\mathbf{R} = \sum_{E\in\mathcal{E}} \mathbf{R}_{E}$.

All software allows the fitting of variance models based on variance
components, for which $\mathbf{G}_U = \sigma^2_U \mathbf{I}_{m(U)}$, 
where $m(U)$ is the number of levels of $U$, and
$\mathbf{R}_E = \sigma^2_E \mathbf{I}_{\Omega_1}$.
The variance-components model for the matrix~$\mathbf{V}$
in equation~\eqref{eq:margmixmodel} is
$\mathbf{V} = \sum_{U \in \mathcal{U}} \sigma^2_U \mathbf{S}_U +
\sum_{E\in\mathcal{E}} \sigma^2_E \mathbf{I}_{\Omega_1}$,
where $\mathbf{Z}_U \mathbf{Z}_U' = \mathbf{S}_U$. For such models,
it is required that all variance components are nonnegative,
which implies that $\mathbf{V}$ is p.s.d.

Some software allows negative estimates of the variance components:
this is essentially fitting a canonical-components model for the
variance matrix, whose general form  is obtained by replacing each
$\sigma^2$-parameter with a $\phi$-parameter:
\begin{equation}
\label{eq:ccmodel}
\mathbf{V} = \sum_{U \in \mathcal{U}} \phi_U \mathbf{S}_U + 
\sum_{E\in\mathcal{E}} \phi_E \mathbf{I}_{\Omega_1}.
\end{equation}
This differs from the variance-components model in that the $\phi_U$,
for $U$ in $\mathcal{U}$,  are not required to be nonnegative,
although $\mathbf{V}$ is required to be p.s.d.

Randomization-based models  are inherently marginal linear mixed models.
The expectation is as given in equation~\eqref{eq:margmixmodel},
with $\mathbf{X}$ replaced by $\mathbf{X}_s\mathbf{X}_h$.
For the variance part of the model,
equation~\eqref{eq:twocancomp} is generalized to
\[
\mathbf{V} = \sum_{i = 1}^p \mathbf{V}_i =
\sum_{i = 1}^p \sum_{H\in\mathcal{H}_i} \phi_H \mathbf{S}_{H},
\]
which is of the form given in equation~\eqref{eq:ccmodel}.

The particular features of randomization-based models are:
\begin{enumerate}
\item[(R.a)] For $i=1$, \ldots, $p$,
all the factors initially defined on $\Omega_i$ are deemed random.
\item[(R.b)] For $i=1$, \ldots, $p$,
$\mathbf{V}_i$ is p.s.d., so that linear combinations of the canonical
  components corresponding to its spectral components must be
  nonnegative. 
This implies that $\phi_{\Omega_i} \geq 0$ for $i = 1$, \ldots, $p$ but
that other canonical components can be negative.
\item[(R.c)] The factors on $\Gamma$, the treatment factors, are
  usually regarded as fixed. 
\end{enumerate}

The set $\mathcal{V}$ of possible estimates of a model for $\mathbf{V}$ is a
subset of the set $\mathcal{M}$ of
$\Omega_1 \times \Omega_1$ real symmetric matrices.
If $\mathcal{F}$ is any subset of $\mathcal{H}$, put
$\mathcal{M}(\mathcal{F}) = \setof{\sum_{F \in \mathcal{F}} a_F \mathbf{S}_F :
a_F \in \mathbb{R} \mbox{ for $F$ in $\mathcal{F}$}}$.
For the models described above, the sets of possible estimates are:
\begin{description}
\item[Unstructured:]  $\mathcal{V}_\mathrm{US}
 = \setof{\mathbf{M} \in \mathcal{M} : \mathbf{M} \mbox{ is p.s.d.}}$;
\item[Canonical-components:]  $\mathcal{V}_\mathrm{CC}(\mathcal{H})
= \mathcal{M}(\mathcal{H}) \cap \mathcal{V}_{\mathrm{US}}$;
\item[Variance-components:]  $\mathcal{V}_\mathrm{VC}(\mathcal{H})
 = \setof{\sum_{H \in \mathcal{H}} a_H \mathbf{S}_H :
a_H \in \mathbb{R}^+_0 \mbox{ for $H$ in $\mathcal{H}$}}$;
\item[Randomization-based:]  
$\mathcal{V}_\mathrm{RB}(\mathcal{H}_1,\ldots,\mathcal{H}_p) =
\setof{\sum_{i = 1}^p \mathbf{V}_i:
\mathbf{V}_i \in \mathcal{V}_{\mathrm{CC}}(\mathcal{H}_i)
\mbox{ for $i=1, \ldots, p$}}$.
\end{description}
Clearly, $\mathcal{V}_\mathrm{VC}(\mathcal{H})
\subset \mathcal{V}_\mathrm{RB}(\mathcal{H}_1,\ldots, \mathcal{H}_p) 
\subset \mathcal{V}_\mathrm{CC}(\mathcal{H}) \subset
\mathcal{V}_\mathrm{US}$.

\subsection{Mixed-model fitting}
\label{sec:mixmod}

By mixed-model fitting, in the case where variances are unknown, we mean
REML estimation of variance parameters followed by EGLS estimation of the
fixed effects. It is preferred for estimation of effects in cases that
are not anova-applicable, including all those without structure 
balance, and for estimation of canonical components when there is not OVS.
It might also be deployed in anova-applicable cases because of software
availability or because it is convenient to use a method that covers
virtually all the cases.
Mixed-model fitting can also be used when $\mathbf{V}$ is known: the
variance parameters are fixed at their known values.
It cannot be used for those anova-applicable cases in which
$\mathcal{R}$ is not orthogonal in relation to $\mathcal{Q}_1$ and
separate analyses are 
required for different parts of $\mathcal{Q}_1$. 
Further advantages of mixed-model fitting are that pseudofactors
are unnecessary and that combined estimates of treatment
effects are obtained when $\mathcal{R}$ is not orthogonal to
the other structures. A disadvantage of mixed-model fitting is that it
is an iterative procedure that can have computational
difficulties. Using anova estimates of canonical components as initial
values helps surmount these. 

In obtaining the fitted values for a randomization-based model using
data from an experiment,  a problem is that mixed-model software
usually fits only variance-components models and perhaps
canonical-components models. The default for GenStat directives
\citep{PayneR12} is to fit canonical-components models and it is an
option in both ASReml-R \citep{Butler09}, a commercial package for R
\citep{R12}, and in \texttt{PROC MIXED} in SAS \citep{SAS10}. The R
packages \texttt{lme4} \citep{Bates12} and \texttt{nlme}
\citep{Pinheiro00} fit   variance-components models only.
Because 
$\mathcal{V}_\mathrm{RB}(\mathcal{H}_1, \ldots, \mathcal{H}_p)
\not\subseteq \mathcal{V}_\mathrm{VC}(\mathcal{H})$,
we recommend fitting canonical-components models.
Even so, a number of difficulties arise:
(i)~all canonical components in the given model must be estimable,
(ii)~software does not allow the separate specification of the
factor sets $\mathcal{H}_i$ and so cannot impose the constraint
that each $\mathbf{V}_i$ is p.s.d.,
and~(iii) software requires
the fitted canonical components to be nonzero 
to avoid singularities in the matrices involved in the computations.

For some variance models, it is inherently impossible to estimate all
variance parameters, as required in~(i).  This should 
be investigated when designing an experiment so that any problems can
be identified and rectified before the experiment is run. The anova
table is extremely useful for this and, in particular, we advocate the
use of `dummy analyses' in which anova is applied to randomly
generated data to check the properties of designs. These analyses are
implemented in GenStat for structure-balanced three-tier experiments;
they can be achieved in R, in many cases including those with
$\mathcal{R}$ not structure balanced in relation to $\mathcal{P}_p$,
by judicious specification of the \texttt{Error} function in a call to
the \texttt{aov} function. Another possibility is to perform a dummy
mixed-model analysis, although this requires that the generated data
are a reasonable fit to the model to be tested. 

Inestimability of variance parameters can arise in two ways. First, a
variance parameter is not estimable if it is 
completely confounded by one or more fixed effects.
For example, if there is some $H$ in $\mathcal{H}$ for which
$\mathbf{S}_H\mathbf{I}_{\mathcal{R}} = \mathbf{S}_H$ then
neither $\phi_H$ nor $\xi_H$ is estimable.
In particular, as pointed out in Section~\ref{sec:inesti}, components
corresponding to the overall mean are never estimable and
mixed-model-fitting software usually excludes them; if not, they must
be dropped.  However, if other canonical components are inestimable 
for this reason then this is usually a sign that the experiment
suffers from some form of pseudoreplication; dropping such components
results in incorrect estimates of standard errors and so is
inadvisable. The avoidance of such design deficiencies is one reason
we urge the use of dummy analyses to check proposed designs. 

Let $\tilde{\mathcal{H}}$ be the set of all factors in $\mathcal{H}$
which do not correspond to the overall mean.
The other cause of inestimability is linear dependence among the matrices
$\mathbf{S}_H$ for $H$ in $\tilde{\mathcal{H}}$. 
Then canonical components need to be dropped so that
those remaining correspond to a linearly independent set of
$\mathbf{S}_H$ for $H$ in $\mathcal{H}^*$,
where $\mathcal{H}^* \subset \tilde{\mathcal{H}}$ 
and $\mathcal{V}_\mathrm{CC}(\mathcal{H}^*) =
\mathcal{V}_\mathrm{CC}(\tilde{\mathcal{H}})$.
The model based on $\mathcal{H}^*$ is a `model of convenience'.
There is a choice about which canonical components to drop in forming
$\mathcal{H}^*$. When all the structures on $\Omega_i$, $i =
1,\dots,p$, exhibit structure balance then a skeleton anova table can
aid in detecting the cause of inestimability of the type outlined in
Section~\ref{sec:inesti} and so in determining which canonical
components to drop. 

Section~\ref{sec:inesti} shows that in Example~\ref{eg:tbb}
the canonical components for Sessions  and Blocks are inestimable.
The term for one or other must be omitted from the mixed model.
This should not be taken to imply that the designer or the analyser of
the experiment is assuming that either does not contribute to the
variability. Indeed, the estimated component should be regarded as
estimating the sum of these two canonical components. Also, the
spectral components for Sessions and Blocks are confounded
and so it is not possible to check that each is nonnegative,
although the nonnegativity of $\xi_\mathrm{S} + 12 \eta_\mathrm{B}$
should be checked. All the other spectral components except $\xi_0$
and $\eta_0$ are estimable and so their nonnegativity can be verified.

For Example~\ref{eg:cotton}, the symbolic
mixed model, derived using Step~1 in \citep{Brien091}, is:
\begin{center}
K $\modsep$ Operatives $\add$ Operatives$\meet$Tests \linebreak $\add$ Blocks $\add$ Blocks$\meet$Plots $\add$ Blocks$\meet$Plots$\meet$Fibres.
\end{center}
As outlined in Section~\ref{sec:inesti}, none of the spectral components
is estimable and so their nonnegativity cannot be checked,
although the nonnegativity of sums such as 
$\xi_{\mathrm{OT}} + \eta_{\mathrm{B}}$ should be.
However, all canonical components other than $\phi_0$ and $\psi_0$
are estimable, except that only $\phi_{\mathrm{OT}} +
\psi_{\mathrm{BPF}}$ is estimable. 
One of these two terms needs to be omitted.
Omitting Blocks$\meet$Plots$\meet$Fibres is, in effect, setting
$\psi_\mathrm{BPF}=0$, and hence $\eta_\mathrm{BPF}=0$.  Again, it is
not assumed that this is the true value of the components. The
constraint is imposed merely to obtain a solution, and
the supposed estimate of $\phi_{\mathrm{OT}}$ is  actually
an estimate of $\phi_{OT} +\psi_{\mathrm{BPF}}$.

For~(ii), a check that the spectral components are
nonnegative is the only option to ensure that the constraints on them are met.
A Gen{S}tat procedure for this has been developed.
Equation~\eqref{eq:zetaphi} is used to obtain the estimated spectral
components from the estimated canonical components. If any 
spectral component is negative then the linear combination of
canonical components on the right-hand side of
equation~\eqref{eq:zetaphi} has to be constrained to zero in a refit
of the model. If there are several negative spectral components, it
may be that some canonical components are constrained to zero.

Difficulty~(iii) occurs because the estimate of some canonical
component happens to be zero. It has to be addressed by removing this
canonical component from the model, which 
implies that a pair of spectral components are equal. It is not
something that can be anticipated ahead of having the data.

\section{Statistical inference}
\label{sec:infer}

In order to perform hypothesis tests or compute confidence intervals
one has to assume that the response 
follows a multi\-variate normal distribution whose expectation and
variance are those described in Section~\ref{sec:multi} 
for the randomization-based model. Some justification for this approach is
that, over all possible randomizations, the distribution of the data
has this expectation and variance. 
The only further assumption that is required
for inference is that of multivariate normality, although the
guarantee for the associated expectation and covariance strictly
applies only over future re-runs of the experiment. The role for
randomization in an analysis based on this model is to ensure that the
sources of variation taken into account by the designer have terms in
the model; that is, it links the model to the design.
Irrespective of the number of tiers, the randomization does not itself
produce distributions whose
third and higher-order moments are those of a multivariate
normal distribution.

\section{Other Models}
\label{sec:altmod}

Steps~2 and~3 of the method in \citep{Brien091} suggest changes that
could be made to the expectation and variance of randomization-based models. 
Here we concentrate on changing treatment factors from fixed to random and 
changing unrandomized factors from random to fixed. The first of these 
produces a randomization-based model, but the second does not; the latter does 
not preserve the variance matrix under randomization as part of the model.

\subsection{Treatment factors regarded as random}
\label{sec:fixtoran}
The simplest modification to the model in equation~\eqref{eq:add} is
to assume that the $\tau_i$, for $i$ in $\Gamma$, are random variables
with common mean $\mu$ and variance matrix $\mathbf{C}_\Gamma$, which
may be as simple as $\sigma^2_\Gamma \mathbf{I}_\Gamma$ or may be
based on a poset block structure on $\Gamma$. So long as $h$~is
equireplicate, $\mathbf{C}_\Gamma$ translates easily to add an extra
variance matrix to $\mathbf{V}$. 

Varieties in early generation variety trials are often
regarded as random: see \cite{Smith01b}.
If Lines are designated random in Example~\ref{eg:Wheat2}
then the variance matrix becomes
\begin{eqnarray*}
\mathbf{V} & = &\xi_0 \mathbf{P}_0 + \xi_\mathrm{O} \mathbf{P}_\mathrm{O}
+ \xi_\mathrm{OIR} \mathbf{P}_\mathrm{OIR}
+ \xi_\mathrm{OIT} \mathbf{P}_\mathrm{OIT}
+ \xi_\mathrm{OIRT} \mathbf{P}_\mathrm{OIRT} \\
 & & + \eta_0 \mathbf{Q}_0 + \eta_\mathrm{B} \mathbf{Q}_\mathrm{B}
+ \eta_\mathrm{BP} \mathbf{Q}_\mathrm{BP}
+ \eta_\mathrm{BPS} \mathbf{Q}_\mathrm{BPS} \\
& & + 8 \sigma^2_\mathrm{L} \mathbf{I}_\mathcal{R}.
\end{eqnarray*}

\subsection{Unrandomized factors regarded as fixed}
\label{sec:rantofix}

Sometimes it is appropriate to classify unrandomized factors
such as Sites, Centres, Laboratories, Sex or Judges as fixed.
It requires that there is no confounding between fixed sources.
It results in the exclusion  of the corresponding subspaces
from the REML estimation of variance parameters,
with canonical components effectively being set to zero and
effects added to the expectation,
so the variance matrix may have LDCVS.
In the expected mean squares,
$q(H)$ replaces $r_ik_H\phi_H$ if generalized factor
$H$ on $\Omega_i$ is designated as fixed.

Suppose that Operatives in Example~\ref{eg:cotton} is to
be considered fixed. This removes $\phi_{\mathrm{O}}$ from
the expression for the variance matrix, and
$\image{(\mathbf{P}_\mathrm{0} + \mathbf{P}_\mathrm{O})}$
is excluded from the REML estimation of the canonical components. The
effect on the expected mean squares in Table~\ref{tab:cox4} is
to replace $\xi_{\mathrm{O}}$ by $\xi_{\mathrm{OT}} + q(\mathrm{O})$ and
$15\phi_\mathrm{O}$ by $q(\mathrm{O})$.

\section{Discussion} \label{sec:discuss}

This paper extends randomization-based models to multitiered experiments with
two or more randomizations in a chain, and discusses the estimation of
treatment effects and their standard errors, and canonical components,
under the assumption of such a model.
There are novel aspects to the estimability of spectral and canonical
components in 
such experiments, including that the variance matrix can exhibit LDCVS.

We have emphasised the usefulness of a skeleton anova in checking the
properties of a design and of anova in analysing anova-applicable
experiments and for supplying initial estimates for mixed-model
fitting. A limitation is software availability. 

Otherwise, mixed-model fitting software is used to
fit a randomization-based model. In this, one has to ensure that
estimates of `variance components' can be negative
and be vigilant that estimates of spectral components are nonnegative.

While potentially negative canonical components are mandated for
randomization-based models, they have the additional benefit
of allowing for negative correlation,
which is realistic in some circumstances: see \citep{JAN:neg}.
\citet{Littel06}, Section~4.7, recommend that unconstrained estimates
be allowed in order to control Type I error, and show that they can
achieve greater power; this agrees with the conclusions of
Wolde-Tsadik and Afifi \citep{WoldeAfifi}. 
However, caution is required in ascribing a negative estimate for a
component to negative population correlation. As Searle, Casella and
McCulloch \citep{Searle92}, Section 3.5, show, for a variance
component just above zero, there can be a high probability of a
negative estimate if the number of treatments is less than 5 and the
number of replicates less than 25. Gilmour and Goos \citep{SGG:poor}
demonstrate that simply allowing negative variance components is not a
panacea, especially in small experiments.

\bibliography{arefs}

\begin{thebibliography}{99}

\bibitem[\protect\citeauthoryear{Alejandro, Bailey, and Cameron}{Alejandro  et~al.}{2003}]{Alejandro03} Alejandro, P.~P., R.~A. Bailey, and P.~J. Cameron (2003). \newblock Association schemes and permutation groups. \newblock {\em Discrete Math.\/}~{\em 266}, 47--67. 
\MR{1991706}

\bibitem[\protect\citeauthoryear{Bailey}{Bailey}{1981}]{Bailey81} Bailey, R.~A. (1981). \newblock A unified approach to design of experiments. \newblock {\em J. R. Statist. Soc. \rm A,\/}~{\em 144}, 214--223. 
\MR{0625801}

\bibitem[\protect\citeauthoryear{Bailey}{Bailey}{1991}]{Bailey91} Bailey, R.~A. (1991). \newblock Strata for randomized experiments. \newblock {\em J. R. Statist. Soc. \rm B. Methodol.\/}~{\em 53}, 27--78. 
\MR{1094275}

\bibitem[\protect\citeauthoryear{Bailey}{Bailey}{2008}]{Bailey08} Bailey, R.~A. (2008). \newblock {\em Design of Comparative Experiments}. \newblock Cambridge: Cambridge University Press. 
\MR{2422352}

\bibitem[\protect\citeauthoryear{Bailey, Ferreira, Ferreira, and Nunes}{Bailey
  et~al.}{2012}]{Bailey12}
Bailey, R.~A., S.~S. Ferreira, D.~Ferreira, and C.~Nunes (2012).
\newblock Estimability of variance components when all model matrices commute.
\newblock submitted for publication.

\bibitem[\protect\citeauthoryear{Bailey, Praeger, Rowley, and Speed}{Bailey  et~al.}{1983}]{BaileyPRS83} Bailey, R.~A., C.~E. Praeger, C.~A. Rowley, and T.~P. Speed (1983). \newblock Generalized wreath products of permutation groups. \newblock {\em Proc. Lond. Math. Soc.\/}~{\em 47}, 69--82. 
\MR{0698928}

\bibitem[\protect\citeauthoryear{Bardin and Aza\"{\i}s}{Bardin and  Aza\"{\i}s}{1990}]{BarAz:rand} Bardin, A. and J.-M. Aza\"{\i}s (1990). \newblock Une hypoth\`ese minimale pour une th\'eorie des plans d'exp\'eriences  randomis\'es. \newblock {\em Rev. Stat. Appl.\/}~{\em 38}, 21--41. 
\MR{1080504}

\bibitem[\protect\citeauthoryear{Bates, Maechler, and Bolker}{Bates
  et~al.}{2012}]{Bates12}
Bates, D., M.~Maechler, and B.~Bolker (2012).
\newblock {\em {lme4}: Linear mixed-effects models using S4 classes}.
\newblock {URL:}
  {\url{http://cran.at.r-project.org/web/packages/lme4/index.html}}, (R package
  version 0.999999-0, accessed November 21, 2012).

\bibitem[\protect\citeauthoryear{Brien}{Brien}{1983}]{Brien83}
Brien, C.~J. (1983).
\newblock Analysis of variance tables based on experimental structure.
\newblock {\em Biometrics,\/}~{\em 39}, 51--59.

\bibitem[\protect\citeauthoryear{Brien}{Brien}{1992}]{Brien92s}
Brien, C.~J. (1992).
\newblock {\em Factorial Linear Model Analysis}.
\newblock Ph.\ D. thesis, Department of Plant Science, The University of
  Adelaide, Adelaide, South Australia.
\newblock {URL:} {\url{ http://hdl.handle.net/2440/37701}}, (accessed November
  21, 2012).

\bibitem[\protect\citeauthoryear{Brien and Bailey}{Brien and  Bailey}{2006}]{BrBa:mult} Brien, C.~J. and R.~A. Bailey (2006). \newblock Multiple randomizations. \newblock {\em J. Roy. Statist. Soc. B Stat. Methodol.\/}~{\em 68}, 571--609. 
\MR{2301010}

\bibitem[\protect\citeauthoryear{Brien and Bailey}{Brien and  Bailey}{2009}]{Brien092} Brien, C.~J. and R.~A. Bailey (2009). \newblock Decomposition tables for experiments. {I.} {A} chain of  randomizations. \newblock {\em Ann. Statist.\/}~{\em 37}, 4184--4213. 
\MR{2572457}

\bibitem[\protect\citeauthoryear{Brien and Bailey}{Brien and  Bailey}{2010}]{Brien10} Brien, C.~J. and R.~A. Bailey (2010). \newblock Decomposition tables for experiments. {II.} {Two-one} randomizations. \newblock {\em Ann. Statist.\/}~{\em 38}, 3164--3190. 
\MR{2722467}

\bibitem[\protect\citeauthoryear{Brien and Dem\'{e}trio}{Brien and  Dem\'{e}trio}{2009}]{Brien091} Brien, C.~J. and C.~G.~B. Dem\'{e}trio (2009). \newblock Formulating mixed models for experiments, including longitudinal  experiments. \newblock {\em J. Agric. Biol. Environ. Stat.\/}~{\em 14}, 253--280. 
\MR{2750840}

\bibitem[\protect\citeauthoryear{Brien, Harch, Correll, and Bailey}{Brien  et~al.}{2011}]{Brien11} Brien, C.~J., B.~D. Harch, R.~L. Correll, and R.~A. Bailey (2011). \newblock Multiphase experiments with at least one later laboratory phase. {I.}  {Orthogonal designs}. \newblock {\em J. Agric. Biol. Environ. Stat.\/}~{\em 16}, 422--450. 
\MR{2843135}

\bibitem[\protect\citeauthoryear{Brien and Payne}{Brien and
  Payne}{1999}]{Brien99}
Brien, C.~J. and R.~W. Payne (1999).
\newblock Tiers, structure formulae and the analysis of complicated
  experiments.
\newblock {\em The Statistician,\/}~{\em 48}, 41--52.

\bibitem[\protect\citeauthoryear{Brien and Payne}{Brien and
  Payne}{2012}]{Brien12AMT}
Brien, C.~J. and R.~W. Payne (2012).
\newblock {AMTIER} {P}rocedure. 
\newblock In 
{\em {GenStat} Reference Manual Release 15, Part 3 Procedures},
  pp.~106--108. Hemel Hempstead, U.K.: VSN International.
\newblock {URL:} \\
  {\url{http://www.vsni.co.uk/resources/documentation/genstat-reference-procedure-library/}},
  (accessed November 21, 2012).

\bibitem[\protect\citeauthoryear{Butler, Cullis, Gilmour, and Gogel}{Butler
  et~al.}{2009}]{Butler09}
Butler, D., B.~R. Cullis, A.~R. Gilmour, and B.~J. Gogel (2009).
\newblock {\em Mixed Models for S language environments: ASReml-R reference
  manual. Version 3}.
\newblock Brisbane: DPI Publications.
\newblock {URL:} {\url{http://www.vsni.co.uk/resources/documentation/}},
  (accessed November 21, 2012).

\bibitem[\protect\citeauthoryear{Cox}{Cox}{1958}]{Cox58b}
Cox, D.~R. (1958).
\newblock The interpretation of the effects of non-additivity in the {Latin}
  square.
\newblock {\em Biometrika\/}~{\em 45}, 69--73.

\bibitem[\protect\citeauthoryear{Curnow}{Curnow}{1959}]{Curnow59}
Curnow, R.~N. (1959).
\newblock The analysis of a two phase experiment.
\newblock {\em Biometrics\/}~{\em 15}, 60--73.

\bibitem[\protect\citeauthoryear{Gilmour and Goos}{Gilmour and  Goos}{2009}]{SGG:poor} Gilmour, S.~G. and P.~Goos (2009). \newblock Analysis of data from non-orthogonal multistratum designs in  industrial experiments. \newblock {\em J. R. Stat. Soc. Ser. C. Appl. Stat.\/}~{\em 58}, 467--484. 
\MR{2750088}

\bibitem[\protect\citeauthoryear{Grundy and Healy}{Grundy and
  Healy}{1950}]{Grundy50}
Grundy, P.~M. and M.~J.~R. Healy (1950).
\newblock Restricted randomization and {Q}uasi-{L}atin squares.
\newblock {\em J. R. Stat. Soc. Ser. B Methodol.\/}~{\em 12}, 286--291.

\bibitem[\protect\citeauthoryear{Houtman and Speed}{Houtman and  Speed}{1983}]{Houtman83} Houtman, A.~M. and T.~P. Speed (1983). \newblock Balance in designed experiments with orthogonal block structure. \newblock {\em Ann. Statist.\/}~{\em 11}, 1069--1085. 
\MR{0720254}

\bibitem[\protect\citeauthoryear{James}{James}{1957}]{James57} James, A.~T. (1957). \newblock The relationship algebra of an experimental design. \newblock {\em Ann. Math. Statist.\/}~{\em 28}, 993--1002. 
\MR{0097142}

\bibitem[\protect\citeauthoryear{Littell, Milliken, Stroup, Wolfinger, and
  Schabenberger}{Littell et~al.}{2006}]{Littel06}
Littell, R.~C., G.~A. Milliken, W.~W. Stroup, R.~D. Wolfinger, and
  O.~Schabenberger (2006).
\newblock {\em SAS for Mixed Models\/} (2nd ed.).
\newblock Cary: SAS Press.

\bibitem[\protect\citeauthoryear{McIntyre}{McIntyre}{1955}]{McIntyre55}
McIntyre, G.~A. (1955).
\newblock Design and analysis of two phase experiments.
\newblock {\em Biometrics\/}~{\em 11}, 324--334.

\bibitem[\protect\citeauthoryear{Monod and Bailey}{Monod and
  Bailey}{1992}]{Monod92}
Monod, H. and R.~A. Bailey (1992).
\newblock Pseudofactors: Normal use to improve design and facilitate analysis.
\newblock {\em J. R. Stat. Soc. Ser. C. Appl. Stat.\/}~{\em 41}, 317--336.

\bibitem[\protect\citeauthoryear{Nelder}{Nelder}{1954}]{JAN:neg} Nelder, J.~A. (1954). \newblock The interpretation of negative components of variance. \newblock {\em Biometrika\/}~{\em 41}, 544--548. 
\MR{0065083}

\bibitem[\protect\citeauthoryear{Nelder}{Nelder}{1965a}]{Nelder65a} Nelder, J.~A. (1965a). \newblock The analysis of randomized experiments with orthogonal block  structure. {I}. {Block} structure and the null analysis of variance. \newblock {\em Proc. Roy. Soc. \textrm{A},\/}~{\em 283}, 147--162. 
\MR{0176576}

\bibitem[\protect\citeauthoryear{Nelder}{Nelder}{1965b}]{Nelder65b} Nelder, J.~A. (1965b). \newblock The analysis of randomized experiments with orthogonal block  structure. {II}. {Treatment} structure and the general analysis of variance. \newblock {\em Proc. Roy. Soc. \textrm{A},\/}~{\em 283}, 163--178. 
\MR{0174156}

\bibitem[\protect\citeauthoryear{Nelder}{Nelder}{1968}]{Nelder68} Nelder, J.~A. (1968). \newblock The combination of information in generally balanced designs. \newblock {\em J. R. Statist. Soc. \rm B. Methodol.\/}~{\em 30}, 303--311. 
\MR{0234582}

\bibitem[\protect\citeauthoryear{Nelder}{Nelder}{1977}]{Nelder77} Nelder, J.~A. (1977). \newblock A reformulation of linear models. \newblock {\em J. R. Statist. Soc. \rm A,\/}~{\em 140}, 48--76. 
\MR{0458743}

\bibitem[\protect\citeauthoryear{Ojima}{Ojima}{1998}]{Ojima98} Ojima, Y. (1998). \newblock General formulae for expectations, variances and covariances of the  mean squares for staggered nested designs. \newblock {\em J. Appl. Stat.\/}~{\em 25}, 785--799. 
\MR{1652264}

\bibitem[\protect\citeauthoryear{Patterson and Thompson}{Patterson and  Thompson}{1971}]{Patterson71} Patterson, H.~D. and R.~Thompson (1971). \newblock Recovery of inter-block information when block sizes are unequal. \newblock {\em Biometrika,\/}~{\em 58}, 545--554. 
\MR{0319325}

\bibitem[\protect\citeauthoryear{Patterson and Thompson}{Patterson and  Thompson}{1975}]{Patterson75} Patterson, H.~D. and R.~Thompson (1975). \newblock Maximum likelihood estimation of components of variance. \newblock In L.~C.~A. Corsten and T.~Postelnicu (Eds.), {\em Proceedings of the  Eighth International Biometric Conference}, pp.\ 197--207. Bucure\c{s}ti:  Editura Academiei Republicii Socialiste Rom\^ania. 
\MR{0468083}

\bibitem[\protect\citeauthoryear{Payne, Welham, and Harding}{Payne
  et~al.}{2012}]{PayneR12}
Payne, R.~W., S.~J. Welham, and S.~A. Harding (2012).
\newblock {\em The Guide to REML in {GenStat}$^{\circledR}$
  (15\textsuperscript{th} Ed.)}.
\newblock Hemel Hempstead: VSN International.
\newblock {URL:} \\ {\url{http://www.genstat.co.uk/resources/documentation/}},
  (accessed November 21, 2012).

\bibitem[\protect\citeauthoryear{Payne and Wilkinson}{Payne and
  Wilkinson}{1977}]{Payne77}
Payne, R.~W. and G.~N. Wilkinson (1977).
\newblock A general algorithm for analysis of variance.
\newblock {\em J. R. Stat. Soc. Ser. C. Appl. Stat.\/}~{\em 26}, 251--260.

\bibitem[\protect\citeauthoryear{Pinheiro and Bates}{Pinheiro and
  Bates}{2000}]{Pinheiro00}
Pinheiro, J.~C. and D.~Bates (2000).
\newblock {\em Mixed Effects Models in S and S-PLUS}.
\newblock New York: Springer.

\bibitem[\protect\citeauthoryear{{R~Development~Core~Team}}{{R~Development~Core~Team}}{2012}]{R12}
{R~Development~Core~Team} (2012).
\newblock {\em R: A Language and Environment for Statistical Computing}.
\newblock Vienna, Austria: R Foundation for Statistical Computing.
\newblock {URL:} \\ {\url{http://www.r-project.org/}}, (accessed November 21,
  2012).

\bibitem[\protect\citeauthoryear{{SAS~Institute~Inc.}}{{SAS~Institute~Inc.}}{2010}]{SAS10}
{SAS~Institute~Inc.} (2010).
\newblock {\em SAS/STAT$^{\circledR}$ 9.22 User�s Guide}.
\newblock Cary, NC: SAS Institute Inc.

\bibitem[\protect\citeauthoryear{Searle, Casella, and Mc{C}ulloch}{Searle  et~al.}{1992}]{Searle92} Searle, S.~R., G.~Casella, and C.~E. Mc{C}ulloch (1992). \newblock {\em Variance {C}omponents}. \newblock New York: John Wiley \& Sons. 
\MR{1190470}

\bibitem[\protect\citeauthoryear{Smith, Cullis, and Gilmour}{Smith  et~al.}{2001}]{Smith01b} Smith, A.~B., B.~R. Cullis, and A.~R. Gilmour (2001). \newblock The analysis of crop variety evaluation data in {A}ustralia. \newblock {\em Aust. N. Z. J. Stat.\/}~{\em 43}, 129--145. 
\MR{1855705}

\bibitem[\protect\citeauthoryear{Speed}{Speed}{1987}]{Speed87a} Speed, T.~P. (1987). \newblock What is an analysis of variance? \newblock {\em Ann. Statist.\/}~{\em 15}, 885--910. 
\MR{0902237}

\bibitem[\protect\citeauthoryear{Speed and Bailey}{Speed and  Bailey}{1987}]{Speed87} Speed, T.~P. and R.~A. Bailey (1987). \newblock Factorial dispersion models. \newblock {\em Int. Stat. Rev.\/}~{\em 55}, 251--277. 
\MR{0963143}

\bibitem[\protect\citeauthoryear{Tjur}{Tjur}{1984}]{Tjur84d} Tjur, T. (1984). \newblock Analysis of variance models in orthogonal designs. \newblock {\em Int. Stat. Rev.\/}~{\em 52}, 33--81. 
\MR{0967202}

\bibitem[\protect\citeauthoryear{Wilkinson}{Wilkinson}{1970}]{Wilkinson70}
Wilkinson, G.~N. (1970).
\newblock A general recursive procedure for analysis of variance.
\newblock {\em Biometrika\/}~{\em 57}, 19--46.

\bibitem[\protect\citeauthoryear{Wolde-Tsadik and Afifi}{Wolde-Tsadik and  Afifi}{1980}]{WoldeAfifi} Wolde-Tsadik, G. and A.~A. Afifi (1980). \newblock A comparison of the ``sometimes pool'', ``sometimes switch'' and  ``never pool'' procedures in the two-way {ANOVA} random effects model. \newblock {\em Technometrics\/}~{\em 22}, 367--373. 
\MR{0585634}

\bibitem[\protect\citeauthoryear{Wood, Williams, and Speed}{Wood
  et~al.}{1988}]{Wood88}
Wood, J.~T., E.~R. Williams, and T.~P. Speed (1988).
\newblock Non-orthogonal block structure in two-phase designs.
\newblock {\em Aust. J. Statist.\/}~{\em 30A}, 225--237.

\bibitem[\protect\citeauthoryear{Yates}{Yates}{1936}]{Yates36}
Yates, F. (1936).
\newblock A new method of arranging variety trials involving a large number of
  varieties.
\newblock {\em J. Agric. Sci.\/}~{\em 26}, 424--455.
\newblock Reprinted with additional author's note in {Yates},~{F}. (1970) {\em
  Experimental} {\em Design}: {\em Selected} {\em Papers}. pages 147--180.
  {Griffin}, {London}.

\end{thebibliography}
\bibliographystyle{chicago}

\end{document}